\documentclass[12pt]{article}

\usepackage{graphicx}%
\usepackage{multirow}%
\usepackage{amsmath,amssymb,amsfonts}%
\usepackage{amsthm}%
\usepackage{mathrsfs}%
\usepackage[title]{appendix}%
\usepackage{xcolor}%
\usepackage{textcomp}%
\usepackage{manyfoot}%
\usepackage{booktabs}%
\usepackage{algorithm}%
\usepackage{algorithmicx}%
\usepackage{algpseudocode}%
\usepackage{listings}%
\usepackage{amsmath}
\usepackage{amsfonts}
\usepackage{amsthm}
\usepackage{graphicx}
\usepackage[utf8]{inputenc}
\usepackage{color}
\usepackage{multirow}
\usepackage{multicol}
\usepackage{booktabs}
\usepackage{blindtext}

\usepackage{graphicx}
\usepackage{tikz}
\usepackage{bm}
\usepackage{float}
\usetikzlibrary{calc}
\usepackage[dvipsnames]{xcolor}

\newcommand{\ma}{\textit{Mathematica$^{\small{\circledR}}$}}

\date{}

\newtheorem{theorem}{Theorem}

\newtheorem{remark}{Remark}

\begin{document}
\title{Differential system related to generalised Krawtchouk polynomials: iterated regularisation and Painlev\'{e} equation}
\author{Galina Filipuk$^a$, Juan Francisco Ma\~{n}as--Ma\~{n}as$^b$,\\ Juan Jos\'{e} Moreno--Balc\'{a}zar$^{b,c}$, Cristina Rodr\'{\i}guez--Perales$^b$*}
\maketitle
{\scriptsize
\noindent $^a$Institute of Mathematics, University of Warsaw, Warsaw, Poland.\\
$^b$Departamento de Matem\'{a}ticas, Universidad de Almer\'{\i}a, Spain.\\
$^c$Instituto Carlos I de F\'{\i}sica Te\'{o}rica y Computacional, Spain.\\
*Corresponding author
\\
Email:\,g.filipuk@uw.edu.pl, jmm939@ual.es, balcazar@ual.es, crp170@ual.es*
}

\begin{abstract}
In this paper we tackle the generalised Krawtchouk polynomials. We show a straightforward connection between certain auxiliary quantities involving the recurrence coefficients of these polynomials and the fifth Painlev\'{e} equation via iterated regularisation. This connection follows an algorithmic procedure avoiding the use of less standardised methods. We further explore how to obtain polynomial systems via iterated regularisation and find decomposition of certain birational transformations into the product of simple ones.
\end{abstract}

\noindent \textbf{Keywords:} Generalised Krawtchouk polynomials $\cdot$  Regularisation  $\cdot$  Painlev\'{e} equations   \\

\noindent \textbf{MSC 2020:}  34M55 $\cdot$ 42C05

\allowdisplaybreaks[4]

\section{Introduction}

It is well known that a sequence of orthogonal polynomials satisfies a three--term recurrence relation (see, for instance,~\cite[p. 18]{Chihara}). In particular, if we consider a sequence of orthonormal polynomials $(p_n)_{n}$, $n=0,1,2,\dots,$ with $p_n(x)=\gamma_n x^n+\dots$, the three--term recurrence relation is given by
$$xp_n(x)=a_{n+1}p_{n+1}(x)+b_n p_n(x)+a_n p_{n-1}(x),$$
where $p_{-1}=0$, $p_0(x)=\gamma_0$, and $a_n$, $b_n$ are known as the recurrence coefficients. Furthermore, the three--term recurrence relation satisfied by the monic orthogonal polynomials $P_n(x)$ given by $P_n(x)=\frac{p_n(x)}{\gamma_n}$  can be expressed as
\begin{equation}
\label{three-term}
xP_n(x)=P_{n+1}(x)+b_nP_n(x)+a_n^2P_{n-1}(x),
\end{equation}
where $P_{-1}=0$, $P_0(x)=1$.  For the classical families of orthogonal polynomials, explicit expressions for their recurrence coefficients are well known. However, for the semiclassical families  such explicit expressions are generally not available. Recently, there has been an increasing interest in studying the connection between the recurrence coefficients of semiclassical orthogonal polynomials and Painlev\'{e} equations (see, for instance~\cite{Van Assche} and the references cited therein).
Frequently, the coefficients $a_n$ and $b_n$ of semiclassical polynomials satisfy nonlinear recurrence relations which, sometimes, can be connected to  Painlev\'{e} equations (see, for instance, \cite{Maxwell Polynomials}, \cite{Boelen1}, \cite{Sextic Freud}--\cite{Truncated Laguerre}, \cite{Class two}, \cite{Min Chen},~\cite{Van Assche2}).  Furthermore, when these coefficients are seen as functions of one of the parameters in the weight, they can be connected to the differential Painlev\'{e} equations.

The famous six Painlev\'{e} equations are nonlinear second--order differential equations which have the Painlev\'{e} property, that is, its general solution
is free from movable branch points. These equations arise in a wide range of mathematical and physical contexts, including integrable systems, statistical physics or random matrix theory (\cite[Chapter 32]{DMLF}). The so--called Painlev\'{e} identification or equivalence problem asks if a given second--order differential equation (possessing the Painlev\'{e} property) can be transformed into one of the six known Painlev\'{e} equations in the canonical form using a change of variables. This problem was solved in \cite{Geom Approach} and \cite{Geom Approach Ham}, in case of birational equivalence, but the identification of the transformation is quite complicated and requires several steps and an extensive use of algebraic geometry. In this paper we show that this procedure can be simplified considerably using the notion of iterated (polynomial) regularisation.  The regularisation procedure and the geometrical aspects of Painlev\'{e} equations have been widely employed to study Painlev\'{e} equations~(see, among others, \cite{Geom Approach,Reg,FR Charlier,KNY, Matano, Okamoto, Sakai}). We also note that our method does not require the full identification procedure employed in the papers \cite{Geom Approach, Geom Approach Ham}. For more information about iterated regularisation we refer to \cite{FR Charlier,Fil Hamil}.

In this paper we focus our attention on a semiclassical generalisation of the Krawtchouk polynomials. Our main objective is to investigate how we can establish connection between some quantities involving the recurrence coefficients of this family of orthogonal polynomials and the fifth Painlev\'{e} equation using iterated regularisation (see Figures~\ref{diagramaIt1} and~\ref{diagramaIt2}). The fifth Painlev\'{e} equation ($P_V$) is given by~\cite[p. 29]{Gromak},~\cite[p. 9]{Van Assche}
\begin{equation}
\label{PV}
y''=\left(\frac{1}{2y}+\frac{1}{y-1}\right)(y')^2-\frac{y'}{t}+\frac{(y-1)^2}{t^2}\left(A_5 y+\frac{B_5}{y}\right)+C_5\frac{y}{t}+D_5\frac{y(y+1)}{y-1},
\end{equation}
where $y=y(t)$ and $A_5,\,B_5,\,C_5,\,D_5$ are arbitrary complex parameters.

For a clearer comprehension of the paper, we first recall the main properties of the classical Krawtchouk polynomials.  The monic Krawtchouk polynomials $P_n(x;\theta,N)$ are given by~\cite[p. 237]{Koekoek book},~\cite[p. 51]{Nikiforov}
\begin{equation*}
P_n(x;\theta,N)=\frac{(-\theta)^n N!}{(N-n)! }\,_2 F_1\left(\begin{array}{c}
-n,-x \\
-N
\end{array} ; \frac{1}{\theta}\right), \quad 0<\theta<1, \ \ N>0,\ \  x\in\{0,\dots, N\},
\end{equation*}
where the hypergeometric function $_rF_s$ is defined by~\cite[p. 5]{Koekoek book}
\begin{equation}
 _r F_s\left(\begin{array}{l}
a_1, \ldots, a_r \\
b_1, \ldots, b_s
\end{array} ;  z\right)=\sum_{k=0}^{\infty} \frac{(a_1)_k(a_2)_k\cdots (a_r)_k}{(b_1)_k(b_2)_k\cdots (b_s)_k}\frac{z^k}{k!},\label{hypergeometric f}
\end{equation}
with the Pochhammer symbol defined by
\begin{equation}
\label{poch}
(a)_0=1, \quad(a)_k=\prod_{i=1}^{k} (a+i-1),\quad k=1,2,3\dots.
\end{equation}
These polynomials are orthogonal with respect to the weight
\begin{equation}
\label{weight kraw}
\omega(x)={N\choose x}\theta^x(1-\theta)^{N-x},
\end{equation}
i.e. they satisfy the following orthogonality condition
$$\sum_{x=0}^N P_n(x;\theta,N) P_m(x;\theta,N){N\choose x}\theta^x(1-\theta)^{N-x}=\frac{N!n!\theta^n(1-\theta)^n}{(N-n)!}\delta_{m,n},$$
where $\delta_{m,n}$ denotes the Kronecker delta. The recurrence coefficients in~(\ref{three-term}) for these polynomials are given by~\cite[p. 120]{Renato Libro},~\cite[p. 44]{Nikiforov},
\begin{equation*}
a_n^2=n\theta(1-\theta)(N-n+1),\quad b_n=N\theta+(1-2\theta)n.
\end{equation*}
\newpage
It is well known that these polynomials correspond to one of the classical discrete families of orthogonal polynomials  (see \cite{Nikiforov}).

A sequence of discrete polynomials with respect to a weight function $\omega(x)$ is said to be classical if there exist polynomials $\sigma(x)$ and $\tau(x)$ with degree $\sigma(x)\leq 2$ and degree $\tau(x)=1$, so that, the following Pearson equation holds
\begin{equation}
\Delta[\sigma(x)\omega(x)]=\tau(x)\omega(x),\label{Pearson}
\end{equation}
where $\Delta$ denotes the forward difference operator $\Delta f(x)=f(x+1)-f(x)$. For the Krawtchouk polynomials we have that~\cite[p. 118]{Renato Libro}
\begin{equation*}
\sigma(x)=x,\quad \tau(x)=\frac{N\theta-x}{1-\theta}.
\end{equation*}

If degree $\sigma(x)> 2$ or degree $\tau(x)\neq1$ in the Pearson equation~(\ref{Pearson}), then the corresponding family of orthogonal polynomials is said to be semiclassical. Obviously the polynomials $\sigma(x)$ and $\tau(x)$ are not unique.  For a family of semiclassical orthogonal polynomials, we define the class $s$ by
$$s=\min_{\sigma,\tau}\max\{\deg \sigma-2, \deg \tau-1\}$$
among all polynomials $\sigma$ and $\tau$ such that the Pearson equation holds (see~\cite{Semiclassical} and the references therein for discrete semiclassical families of orthogonal polynomials).\\

In this paper, we consider a generalisation of the monic Krawtchouk polynomials $P_n(x):=P_n(x;t,\alpha,N)$, which were extensively studied in~\cite{Krawtchouk Polynomials}, orthogonal with respect to the weight
\begin{equation}
\label{weight}
\omega(x):=\omega(x;t,\alpha,N)={N\choose x}\frac{t^x}{(1-\alpha)_x},\quad t>0,\quad\alpha<1,\quad  x\in\{0,\dots, N\},
\end{equation}
where the Pochhammer symbol $(a)_k$ is defined by (\ref{poch}). Observe that for this weight the following Pearson equation holds
\begin{equation*}
\Delta[\sigma(x)\omega(x)]=\tau(x)\omega(x),
\end{equation*}
with
\begin{equation*}
\sigma(x)=(x-\alpha)x,\quad \tau(x)=-x^2-(t-\alpha)x+tN,
\end{equation*}
since
\newpage
\begin{align*}
&\Delta\left[(x-\alpha)x{N\choose x}\frac{t^x}{(1-\alpha)_x}\right]\\
=&(x+1-\alpha)(x+1){N\choose x+1}\frac{t^{x+1}}{(1-\alpha)_{x+1}}-(x-\alpha)x{N\choose x}\frac{t^x}{(1-\alpha)_x}\\
=&\left[(N-x)t-(x-\alpha)x\right]{N\choose x}\frac{t^x}{(1-\alpha)_x}.
\end{align*}

Therefore, this family of polynomials is semiclassical of class $1$. Furthermore, we highlight that when $t\to\infty$, $-\alpha\to\infty$ and $-\frac{t}{\alpha}\to\frac{\theta}{1-\theta}$ in the weight definition~(\ref{weight}) we recover the Krawtchouk weight~(\ref{weight kraw}).

For these generalised monic Krawtchouk polynomials we rewrite the three--term recurrence relation~(\ref{three-term}) as
\begin{equation}
\label{three-term generalised}
xP_n(x)=P_{n+1}(x)+b_n(t)P_n(x)+a_n^2(t)P_{n-1}(x),
\end{equation}
highlighting the dependence on $t$ of the coefficients.

%It can be shown that the quantities $q$ and $p$ involving the recurrence coefficients $a_n^2(t)$ and $b_n(t)$ of the generalised Krawtchouk polynomials
% $$q:=\frac{1}{N}\left(\frac{a_n^2}{t}+n\right),\quad p:=-\frac{b_n+N+1+t-n-\alpha}{N}$$
%satisfy the following differential system (see~\cite{Krawtchouk Polynomials}):
%\begin{equation}
%\label{original system}
%\begin{gathered} %not t,
%q'=\frac{N p^2 (n-N q)+2 N p q (n-N q)+q (N q (\alpha +n-2 N-2)+(N+1) (-\alpha +N+1))}{N t (p+q)},\\
%p'=\frac{2 N p q t+p ((N+1) (-\alpha +N+1)+N p (-\alpha +N p+2 N+t+2))+N q^2 t}{N t (p+q)},
%\end{gathered}
%\end{equation} where $q=q(t),\, p=p(t)$ and $'=\frac{d}{dt}$.

It is known that certain quantities involving these coefficients $a_n^2(t)$ and $b_n(t)$ satisfy a differential system (see~\cite{Krawtchouk Polynomials}). To be self--contained, in Section \ref{obtaining the system} we will detail how to construct this system, which requires making use of a discrete system obtained using the technique of ladder operators \cite{Krawtchouk Polynomials} and the Toda system~\cite[Theorem 3.8]{Van Assche}.

In~\cite{Krawtchouk Polynomials} the authors identified a connection between that differential system and the fifth Painlev\'{e} equation involving too complicated expressions.  In fact, they did not provide the explicit form of the final second--order differential equation, see~\cite[Eq. (40)]{Krawtchouk Polynomials}, due to its complexity. Moreover, their method to deduce the transformation which establishes the connection with the fifth Painlev\'{e} equation relied on guessing some form of the birational transformation, based on analogies with other semiclassical polynomials, especially generalised Meixner polynomials, and then determine computationally the unknown coefficients. However, our contribution consists in the use of a different procedure, the so--called iterated regularisation (see Figures~\ref{diagramaIt1} and~\ref{diagramaIt2}). This method allows us to obtain a connection between that system and the fifth Painlev\'{e} equation in a more systematic way. Thus, the main advantage lies in having an algorithmic procedure, as illustrated in the flow diagram of Figure~\ref{diagrama}. That is, after a sequence of systematic  transformations, we obtain simple final systems which can be identified with $P_V$ in a straightforward way or at most using a M\"{o}bius transformation. In particular, as mentioned previously, whereas in~\cite{Boelen1} the authors obtained a second--order differential without providing its explicit form, in our work we obtain equations (\ref{equationU11})--(\ref{equationU31}), (\ref{eq12})--(\ref{eq31}) which are, or related to, the fifth Painlev\'{e} equation given by~(\ref{PV}). Therefore, this procedure can be implemented in a computational language.

Indeed, we give two alternatives to establish a connection with the fifth Painlev\'{e} equation, the first one only involving the regularisation of the differential system and a second alternative thoughtout polynomial regularisation (see Figure~\ref{diagrama} for more details).

The structure of this paper is as follows. In Section~\ref{obtaining the system}, following~\cite{Krawtchouk Polynomials}, we explain how system~\eqref{original system2} can be deduced. Section~\ref{regularisation} is devoted to the regularisation of this system and iterated regularisation is studied in Section~\ref{iterated regularisation}. In Section~\ref{conection PV}, we address how iterated regularisation enables us to study the connection to $P_V$. In Section~\ref{polynomial reg}, we apply one more iteration of the regularisation procedure and obtain polynomial systems which are reduced directly to the fifth Painlev\'{e} equation.  Here, we also include a colored graph (Figure~\ref{diagrama}) with the main information of Sections \ref{regularisation} and \ref{iterated regularisation}, summarising the regularisation process and the connections with the fifth Painlev\'{e} equation. Finally, in Section~\ref{hamiltonian} we observe that the regular systems obtained in previous sections are Hamiltonian. All the computations presented in this work were carried out using \ma  \textit{14.1.}

\section{The differential system}
\label{obtaining the system}
In this section, in order to be self--contained,  we recall the procedure to obtain the system~\eqref{original system2}. This system results from combining a discrete system for the recurrence coefficients in~\eqref{three-term generalised} for generalised Krawtchouk polynomials~\cite[Theorem 1.1]{Krawtchouk Polynomials} and the Toda--type differential discrete system~\cite[Theorem 3.8]{Van Assche}. These results are given in the following statements. \\

\begin{theorem}~\cite[Theorem 1.1]{Krawtchouk Polynomials}
Let $a_n^2$  and $b_n$ be the recurrence coefficients in~\eqref{three-term generalised} for the weight~\eqref{weight}. The following quantities
\begin{equation}
\label{quantities xn yn}
x_n:=\frac{1}{N}\left(\frac{a_n^2}{t}+n\right),\quad y_n:=-\frac{b_n+N+1+t-n-\alpha}{N},
\end{equation}
satisfy the discrete system
\begin{equation}
\label{discrete system}
\left\{
  \begin{array}{ll}
    (x_n+y_n)(x_{n+1}+y_{n})=\displaystyle-\frac{y_n(N+1+Ny_n)(N+1-\alpha+Ny_n)}{tN} ,\\ \\
    (x_n+y_n)(x_n+y_{n-1})=\displaystyle\frac{x_n(-N-1+Nx_n)(\alpha-N-1+Nx_n)}{N(Nx_n-n)},
  \end{array}
\right.
\end{equation}
with initial conditions
$$x_0=0,\quad y_0=-\frac{N+1+t-\alpha}{N}-\frac{t}{1-\alpha}\frac{M(-N+1,2-\alpha,-t)}{M(-N,1-\alpha,-t)},$$
where $M(a,b,z)$ is the confluent hypergeometric function $_1F_1(a;b;z)$ defined by~(\ref{hypergeometric f}).\\
\end{theorem}

\begin{theorem}~\cite[Theorem 3.8]{Van Assche}
\label{toda}
Suppose $\mu$ is a positive measure on the real line for which all the moments exist, and $\mu_z$ is the modified measure with $d{\mu}_z(x)=e^{xz}d\mu(x),$ where $z\in\mathbb{R}$ is such that all the moments of $\mu_z$ exist. Then the recurrence coefficients $a_n^2$ and $b_n$ of the orthogonal polynomials for the measure $\mu_z$ satisfy
\begin{align*}
\frac{d}{dz}a_n^2&=a_n^2(b_n-b_{n-1}),\quad n\geq 1,\\
\frac{d}{dz}b_n&=a_{n+1}^2-a_n^2,\quad n\geq 0,
\end{align*}
with $a_0^2=0.$
\end{theorem}

To deduce system~\eqref{original system2}, we first observe that setting $t=e^z$ in the weight~\eqref{weight} we have
$$\omega(x)={N\choose x}\frac{t^x}{(1-\alpha)_x}=e^{zx}\frac{1}{(1-\alpha)_x}{N\choose x}.$$
Thus, according to Theorem~\ref{toda} the recurrence coefficients of the orthogonal polynomials with respect to the weight~\eqref{weight} satisfy the following Toda--type system
\begin{align*}
\frac{d}{dt}a_n^2&=\frac{a_n^2}{t}(b_n-b_{n-1}),\quad n\geq 1,\\
\frac{d}{dt}b_n&=\frac{1}{t}(a_{n+1}^2-a_n^2),\quad n\geq 0,
\end{align*}
where we have used $dt=e^z dz=tdz.$ Moreover, we shall express this Toda--type system in terms of the quantities $x_n, y_n$ defined by~\eqref{quantities xn yn}. Taking into account that $a_n^2=(Nx_n-n)t$ and $b_n=-Ny_n-N-1-t+n+\alpha,$ we obtain
\begin{align*}
\frac{d}{dt}\left((Nx_n-n)t\right)&=(Nx_n-n)(-Ny_n+Ny_{n-1}+1),\\
\frac{d}{dt}(-Ny_n-N-1-t+n+\alpha)&=\frac{1}{t}((Nx_{n+1}-n-1)t-(Nx_n-n)t),
\end{align*}
which is equivalent to
\begin{equation}
\begin{gathered}
\label{toda system xn}
\frac{d}{dt}x_n=\frac{(Nx_n-n)(-Ny_n+Ny_{n-1})}{tN},\\
\frac{d}{dt}y_n=x_n-x_{n+1}.
\end{gathered}
\end{equation}
Now, we find the following values for $x_{n+1}$ and $y_{n-1}$ from the first and second equations of~\eqref{discrete system}, respectively,
\newpage
\begin{align*}
x_{n+1}&=-\frac{y_{n}(N+1+Ny_n)(N+1-\alpha+Ny_n)+y_n(x_n+y_n)tN}{(x_n+y_n)tN},\\
y_{n-1}&=\frac{x_n(-N-1+Nx_n)(\alpha-N-1+Nx_n)-x_nN(Nx_n-n)(x_n+y_n)}{N(Nx_n-n)(x_n+y_n)}.
\end{align*}
Substituting these values into~\eqref{toda system xn} we obtain
\begin{align}
\frac{d}{dt}x_n&=\frac{(Nx_n-n)}{tN}\left(-Ny_n+\frac{x_n(-N-1+Nx_n)(\alpha-N-1+Nx_n)-x_nN}{(Nx_n-n)(x_n+y_n)}\right)\notag\\
&=\frac{(Nx_n-n)}{t}\left(-y_n-x_n+\frac{x_n(-N-1+Nx_n)(\alpha-N-1+Nx_n)}{N(Nx_n-n)(x_n+y_n)}\right)\label{toda1xn},
\end{align}
and
\begin{equation}
\frac{d}{dt}y_n=\, x_n+y_n+\frac{y_n(N+1+Ny_n)(N+1-\alpha+Ny_n)}{tN(x_n+y_n)}\label{toda1yn}.
\end{equation}
Finally, from equations~\eqref{toda1xn} and~\eqref{toda1yn}, renaming the variables as $p(t):=x_n$, $q(t):=y_n$ for an arbitrary fixed $n$ and simplifying we obtain the system

\begin{equation}
\label{original system2}
\begin{gathered} %not t,
q'=\frac{N p^2 (n-N q)+2 N p q (n-N q)+q (N q (\alpha +n-2 N-2)+(N+1) (-\alpha +N+1))}{N t (p+q)},\\
p'=\frac{2 N p q t+p ((N+1) (-\alpha +N+1)+N p (-\alpha +N p+2 N+t+2))+N q^2 t}{N t (p+q)},
\end{gathered}
\end{equation}
where $'$ means $\frac{d}{dt}$.

\section{Regularisation of the system}
\label{regularisation}
In this section we address the regularisation of system~\eqref{original system2}. To begin we provide an overview of this procedure. We consider the system in the $(q,p)$ coordinates on the complex projective space $\mathbb{P}^1\times \mathbb{P}^1$, since it is necessary to take into account the poles of $q(t)$ and $p(t)$. Therefore, we also analyse the system in the coordinate charts $(q,P)$, $(Q,p)$ and $(Q,P),$ where $Q={1}/{q},\,P={1}/{p}$.

The regularisation process consists in identifying the points of indeterminacy of  system~\eqref{original system2} within these four charts and resolving them using the blow--up procedure (see \cite{Shafarevich} for more details). Points of indeterminacy arise when both the numerator and denominator on the right--hand side of at least one equation vanish, leading to the vector field undefined. To blow up a point $P_*=(a,b)$, two new charts $(u,v),\, (U, V)$ are introduced via the following changes of variable $x=a+uv=a+V,\,y=b+v=b+UV.$ In these charts, the exceptional divisors are given by $v=0$ and $V=0$.

By resolving a sequence of points of indeterminacy, we finally obtain regular systems on the final exceptional divisors, which we refer to as final chart systems. The steps of the regularisation procedure are graphically described in Figure~\ref{diagramaIt1}.

%%%%%%%%%%%%%%%%%%%%%%%%%%%%%%%%%%%%%%%%%%%%%%%%%%%%%%%%%%%%%%%%%%%%%%%%%%%%%%%%%%%%%% Diagram of regularisation algorithm

\begin{figure}[!h]
\hspace{7em}
{\footnotesize
\begin{tikzpicture}[>=latex]
\label{flowchartit1}
  %
  % Styles for states, and state edges
  %
  \tikzstyle{state} = [draw, very thick, fill=white, rectangle, minimum height=2em, minimum width=2em, node distance=5em, font={\sffamily}]
  \tikzstyle{stateit1} = [draw, very thick, fill=white, rectangle, minimum height=2em, minimum width=2em, node distance=6em, font={\sffamily}]
  \tikzstyle{stateEdgePortion} = [black,very thick];
  \tikzstyle{stateEdge} = [stateEdgePortion,->];  % sin flechas
%  \tikzstyle{stateEdge} = [stateEdgePortion,->];  % con flechas
  \tikzstyle{edgeLabel} = [pos=0.5, text centered,very thick, font={\sffamily\small}];

\hspace{-2cm}\node[state, name=a1] {\textbf{Initial system on $\mathbb{P}^1\times \mathbb{P}^1$}};

\node[state, name=a2,below of=a1] {Points of indeterminacy};
\draw ($(a1.south) + (0em,0)$)
      edge[stateEdge] node[edgeLabel]{}
      ($(a2.north)$);

\node[state, name=a2b,below of=a2] {Blow up of the point};
\draw ($(a2.south) + (0em,0)$)
      edge[stateEdge] node[edgeLabel]{}
      ($(a2b.north)$);

\node[state, name=a3,below of=a2b] {Check the system regularity};
\draw ($(a2b.south) + (0em,0)$)
      edge[stateEdge] node[edgeLabel]{}
      ($(a3.north)$);

\node[stateit1, name=a4,below of=a3,left of=a3, xshift=-5em] {The regularisation algorithm stops};
\draw ($(a3.south) + (0em,0)$)
      edge[stateEdge] node[edgeLabel,xshift=-2em,yshift=0em]{Yes}
      ($(a4.north)$);

\node[stateit1, name=a5,below of=a4,right of=a3, xshift=3em] {Repeats the procedure. Find cascades for such points};
\draw ($(a3.south) + (0em,0)$)
      edge[stateEdge] node[edgeLabel,xshift=2em,yshift=0em]{No}
      ($(a5.north)$);
  \draw ($(a5.east) $)
      edge[stateEdge, bend right=75] node[edgeLabel]{}
      ($(a2.east)$);

\end{tikzpicture}
}
\caption{Diagram of the regularisation algorithm.}
\label{diagramaIt1}
\end{figure}
%%%%%%%%%%%%%%%%%%%%%%%%%%%%%%%%%%%%%%%%%%%%%%%%%%%%%%%%%%%%%%%%%%%%%%%%%%%%%%%%%%%%%

Now, we apply this regularisation process to system~\eqref{original system2}. We first compute the following five points of indeterminacy in different charts:
\begin{align}
P_1&=(q=0,\,p=0),\quad P_2=\left(q=\frac{1+N}{N},\,p=-\frac{1+N}{N}\right),\notag\\
P_3&=\left(q=\frac{1+N-\alpha}{N},\,p=-\frac{1+N-\alpha}{N}\right),\quad P_4=\left(q=\frac{n}{N},\,P=0\right),\notag\\
P_5&=(Q=0,\,P=0).\label{points}
\end{align}
We note that points $P_2$ and $P_3$ coincide when $\alpha=0$ (further details on this specific parameter value will be provided later). Blowing up points $P_1$, $P_2$, $P_3$ and $P_4$ we obtain regular systems on the corresponding exceptional divisors. However, after blowing up point $P_5$ the resulting rational system is not regular on the exceptional divisor, requiring the resolution of a sequence of points. We refer the reader to the diagram in Figure~\ref{diagrama} to see a complete graphical representation of the regularisation procedure applied to the initial system~(\ref{original system2}).

To clarify the notation we denote different charts as follows: after blowing up points $P_i$, $i=1,\dots,5$ given by~\eqref{points} we introduce  $(u_{i,1},v_{i,1}),(U_{i,1},V_{i,1})$ charts. Analogously, blowing up a point of indeterminacy of the system in $(u_{i,1},v_{i,1})$, $(U_{i,1},V_{i,1})$ charts, we introduce $(u_{i,2},v_{i,2}),(U_{i,2},V_{i,2})$ charts. %%No abuse of notation may arise since we will only consider the system in one of the two possible charts.
In general, we introduce $(u_{i,j},v_{i,j})$, $(U_{i,j},V_{i,j})$ charts after blowing up a point of indeterminacy of systems in $(u_{i,j-1},v_{i,j-1})$, $(U_{i,j-1},V_{i,j-1})$ charts, for $j\geq 2$. To denote blow--up transformations we will use the following notation, we denote by $\phi_{i,j}$ and $\hat{\phi}_{i,j}$ the blow--up of a point introducing $(u_{i,j},v_{i,j})$ and $(U_{i,j},V_{i,j})$ charts, respectively. For instance, to blow up a point of indeterminacy $P_i=(x_*,y_*)$ of  system~\eqref{original system2} we apply the following transformations
$$\phi_{i,1}:\ q= x_*+u_{i,1} v_{i,1},\quad p= y_*+v_{i,1},$$
$$\hat{\phi}_{i,1}:q=x_*+ V_{i,1},\quad p= y_*+ U_{i,1}V_{i,1},$$
where $q,p$ will change depending on the initial chart in which we are considering system~\eqref{original system2}.

When we resolve a cascade of points, we might need an intermediate change of variables introduced by Takano et al in~\cite{Matano} which we call a twist and denote by $\tau_{i,j}$ or $\hat{\tau}_{i,j}$. This transformation is a local change of variable given by $\tau_{i,j}: u_{i,j}\to \frac{1}{u_{i,j}},\,\,v_{i,j}\to {v_{i,j}}$ or $\hat{\tau}_{i,j}: U_{i,j}\to \frac{1}{U_{i,j}},\,\,V_{i,j}\to {V_{i,j}}$. In the following, we proceed to detail the resolution of each point of indeterminacy.

Blowing up  $P_1=(q=0,\,p=0)$, which involves the following transformations in both charts
$$\phi_{1,1}:\ q= u_{1,1} v_{1,1},\quad p= v_{1,1},$$
\begin{equation}
\label{phi11}
\hat{\phi}_{1,1}:q= V_{1,1},\quad p=  U_{1,1}V_{1,1},
\end{equation}
we obtain the following systems (we remind that $'=\frac{d}{dt}$)
\begin{equation}
\label{systuv11original}
\begin{gathered}
u_{1,1}'=\frac{n-u_{1,1}^2t-\left(u_{1,1} \left(-\alpha -n+2 N v_{1,1}+2 N+t+2\right)\right)}{t},
\\
v_{1,1}'=\frac{N^2 v_{1,1}^2+(N+1)
   (-\alpha +N+1)+N v_{1,1} \left(-\alpha +2 N+t u_{1,1}^2+2 t u_{1,1}+t+2\right)}{N t \left(u_{1,1}+1\right)},
\end{gathered}
\end{equation}

\begin{equation}
\label{systUV11original}
\begin{gathered}
U_{1,1}'=\frac{U_{1,1} (-\alpha -n+2 N+t+2)-U_{1,1}^2 \left(n-2 N V_{1,1}\right)+t}{t},\\
V_{1,1}'=\frac{N V_{1,1} \left(U_{1,1} \left(U_{1,1}+2\right) \left(n-N V_{1,1}\right)+\alpha +n-2 N-2\right)+(N+1) (-\alpha +N+1)}{N t \left(U_{1,1}+1\right)},
\end{gathered}
\end{equation}
which are regular on $v_{1,1}=0$ and $V_{1,1}=0$, respectively.\\

Similarly, blowing up $P_2=\left(q=\frac{1+N}{N},\,p=-\frac{1+N}{N}\right)$, that is, doing the transformations
$$\phi_{2,1}:\ q= \frac{1+N}{N}+u_{2,1} v_{2,1},\quad p= -\frac{1+N}{N}+ v_{2,1},$$
\begin{equation}
\label{phi21}
\hat{\phi}_{2,1}:\ q= \frac{1+N}{N}+V_{2,1},\quad p= -\frac{1+N}{N}+ U_{2,1}V_{2,1},
\end{equation}
the following two systems in different charts are obtained
\begin{equation}
\label{systuv21original}
\begin{gathered}
u_{2,1}'=\frac{u_{2,1} \left(\alpha +n-2 N v_{2,1}-t\right)+n-N-t u_{2,1}^2-1}{t},\\
v_{2,1}'=\frac{N^2 v_{2,1}^2+\alpha  (N+1)-N v_{2,1}
   \left(\alpha +N-t u_{2,1}^2-2 t u_{2,1}-t+1\right)}{N t \left(u_{2,1}+1\right)},
\end{gathered}
\end{equation}

\begin{equation}
\label{systUV21original}
\begin{gathered}
U_{2,1}'=\frac{U_{2,1}^2 \left(-n+2 N V_{2,1}+N+1\right)-U_{2,1} (\alpha +n-t)+t}{t},\\
V_{2,1}'=\frac{N V_{2,1} \left(-U_{2,1} \left(U_{2,1}+2\right) \left(N V_{2,1}-n+N+1\right)+\alpha +n\right)+\alpha  (N+1)}{N t \left(U_{2,1}+1\right)},
\end{gathered}
\end{equation}
which are regular on $v_{2,1}=0$ and $V_{2,1}=0$ (while $\alpha\neq 0$), respectively.\\
\begin{remark} When $\alpha=0$ systems~\eqref{systuv21original} and~\eqref{systUV21original} are not regular, since they have the following point of indeterminacy $\tilde{P}_{2,2}=(u_{2,1}=-1,v_{2,1}=0)=(U_{2,1}=-1,V_{2,1}=0)$ on the exceptional divisor ${v}_{2,1}=0$ and ${V}_{2,1}=0$, respectively. Thus, to obtain regular systems we need to perform one more blow--up. We resolve this point applying the following transformation in system~\eqref{systuv21original} with $\alpha=0$,
$$u_{2,1}=-1+\tilde{V}_{2,2},\quad v_{2,1}=\tilde{U}_{2,2}\tilde{V}_{2,2},$$
where we use the notation $(\tilde{U}_{2,2},\tilde{V}_{2,2})$ to distinguish that we are working with the specific value of the parameter $\alpha=0$. We obtain the following regular system on the exceptional divisor $\tilde{V}_{2,2}=0$ (we only show the system in one chart to avoid including too many systems):
\begin{equation}
\label{systUV220}
\begin{gathered}
\tilde{U}_{2,2}'=\frac{-n \tilde{U}_{2,2}+2 N \tilde{U}_{2,2}^2 \tilde{V}_{2,2}-N \tilde{U}_{2,2}^2+2 t \tilde{U}_{2,2} \tilde{V}_{2,2}-t \tilde{U}_{2,2}}{t},\\
\tilde{V}_{2,2}'=\frac{n \tilde{V}_{2,2}-2 N \tilde{U}_{2,2} \tilde{V}_{2,2}^2+2 N \tilde{U}_{2,2} \tilde{V}_{2,2}-t \tilde{V}_{2,2}^2+t \tilde{V}_{2,2}-N-1}{t}.
\end{gathered}
\end{equation}

\end{remark}

Now, we proceed in a similar way for point $P_3=\left(q=\frac{1+N-\alpha}{N},\,p=-\frac{1+N-\alpha}{N}\right)$. Blowing up the point, that is, applying the transformations
 $$\phi_{3,1}:\ q= \frac{1+N-\alpha}{N}+u_{3,1} v_{3,1},\quad p= -\frac{1+N-\alpha}{N}+ v_{3,1},$$
 \begin{equation}
\label{phi31}
\hat{\phi}_{3,1}:\ q= \frac{1+N-\alpha}{N}+V_{3,1},\quad p= -\frac{1+N-\alpha}{N}+ U_{3,1}V_{3,1},
\end{equation}
the resulting systems are given by
\begin{equation*}
\begin{gathered}
u_{3,1}'=\frac{ -u_{3,1} \left(\alpha -n+2 N v_{3,1}+t+tu_{3,1}\right)+\alpha+n-N-1}{t},\\
v_{3,1}'=\frac{N^2 v_{3,1}^2+\alpha
   (\alpha -N-1)+N v_{3,1} \left(2 \alpha -N+t u_{3,1}^2+2 t u_{3,1}+t-1\right)}{N t \left(u_{3,1}+1\right)},
\end{gathered}
\end{equation*}

\begin{equation}
\label{systUV31original}
\begin{gathered}
U_{3,1}'=\frac{U_{3,1}^2 \left(-\alpha -n+2 N V_{3,1}+N+1\right)+U_{3,1} (\alpha -n+t)+t}{t},\\
V_{3,1}'=\frac{\alpha  (\alpha -N-1)-N V_{3,1} \left(\alpha +U_{3,1} \left(U_{3,1}+2\right) \left(-\alpha -n+N V_{3,1}+N+1\right)-n\right)}{N t \left(U_{3,1}+1\right)},
\end{gathered}
\end{equation}
which are regular on $v_{3,1}=0$ and $V_{3,1}=0$, respectively. To reduce the amount of systems included, from this point forward, we will only present the resulting systems in just one of the charts.

Blowing up $P_4=\left(q=\frac{n}{N},\,P=0\right)$, we apply the following transformation
\begin{equation}
\label{phi41}
\hat{\phi}_{4,1}:\ q= \frac{n}{N}+V_{4,1},\quad P= U_{4,1}V_{4,1},
\end{equation}
and we obtain the following regular system on $V_{4,1}=0$
\begin{equation}
\label{systUV41original}
\begin{gathered}
U_{4,1}'=-\frac{R_1^{4,1}(U_{4,1},V_{4,1},t)}{N t \left(U_{4,1} V_{4,1} \left(n+N
   V_{4,1}\right)+N\right)},\\
V_{4,1}'=\frac{R_2^{4,1}(U_{4,1},V_{4,1},t)}{N t U_{4,1} \left(U_{4,1} V_{4,1} \left(n+N V_{4,1}\right)+N\right)},
\end{gathered}
\end{equation}
where \\

$R_1^{4,1}(U_{4,1},V_{4,1},t)=U_{4,1} (U_{4,1} (N V_{4,1} (2 (-\alpha -N (\alpha +2 n-2)+n (\alpha
   +n+t-2)+N^2+1)+N V_{4,1} (\alpha +n-2 N+2 t-2))+n (-n+N+1) (-\alpha -n+N+1))+N^2 (-\alpha -2 n-2 N
   V_{4,1}+2 N+t+2)+t U_{4,1}^2 V_{4,1}^2 (n+N V_{4,1}){}^2),$\\

$R_2^{4,1}(U_{4,1},V_{4,1},t)=U_{4,1} V_{4,1} (n+N V_{4,1}) (U_{4,1} ((-n+N+1) (-\alpha -n+N+1)+N V_{4,1}
   (\alpha +n-2 N-2))-2 N^2)-N^3$.\\

In the case of point $P_5=(Q=0,\,P=0)$ we need to perform a sequence of birational transformations to obtain a regular system, where we first introduce a twist. Applying the following sequence of transformations
\begin{align}
\phi_{5,1}&:\, Q= u_{5,1} v_{5,1},\quad P= v_{5,1} \label{phi51},\\
\hat{\phi}_{5,2}&:\, u_{5,1}= V_{5,2},\quad v_{5,1}= U_{5,2}V_{5,2}\label{phi52},\\
\hat{\tau}_{5,2}&:\, U_{5,2}\to \frac{1}{U_{5,2}},\label{tau52}\\
\phi_{5,3}&:\, u_{5,2}= -\frac{N}{t}+u_{5,3} v_{5,3},\quad v_{5,2}= v_{5,3}\label{phi53},\\
\phi_{5,4}&:\, u_{5,3}= \frac{-1-2N+n-t+\alpha}{t}+u_{5,4} v_{5,4},\quad v_{5,3}= v_{5,4}\label{phi54},
\end{align}
we obtain for this cascade a final chart system which is regular on the exceptional divisor $v_{5,4}=0$
\begin{equation}
\label{systuv54op1}
\begin{gathered}
u_{5,4}'=\frac{(N-n) (-\alpha -n+N)-N u_{5,4} \left(-v_{5,4} \left(v_{5,4}+2\right) \left(n-N u_{5,4}\right)+\alpha +n-2 N\right)}{N t
   \left(v_{5,4}+1\right)},\\
v_{5,4}'=\frac{-n v_{5,4}^2+n v_{5,4}+2 N u_{5,4} v_{5,4}^2-2 N v_{5,4}-t v_{5,4}+\alpha  v_{5,4}-t}{t}.
\end{gathered}
\end{equation}

\noindent By composing the previous transformations we deduce that the regularising transformation from the original system~\eqref{original system2} in $(Q,P)$ coordinates to system~\eqref{systuv54op1} is given by
$$\Phi_{5,4}=\phi_{5,4}\circ\phi_{5,3}\circ\hat{\tau}_{5,2}\circ\hat{\phi}_{5,2}\circ\phi_{5,1},$$
\begin{align*}
\Phi_{5,4}:\ Q&= \frac{N v_{5,4}^2}{v_{5,4} \left(N u_{5,4} v_{5,4}+\alpha +n-2 N-t-1\right)-t},\\
 P&=\frac{N v_{5,4}}{v_{5,4} \left(N u_{5,4}
   v_{5,4}+\alpha +n-2 N-t-1\right)-t},
\end{align*}
and its inverse is given by
$$ \Phi_{5,4}^{-1}:\, u_{5,4}= \frac{Q (P (-\alpha -n+2 N+t+1)+N)+t P^2}{N Q^2},\quad v_{5,4}= \frac{Q}{P}.$$

On the other hand, we could also have applied the following sequence of transformations to obtain a regular system
\begin{align}
\phi_{5,1^b}&:\, Q= u_{5,1^b} v_{5,1^b},\quad P= v_{5,1^b},\label{phi51b}\\
\phi_{5,2^b}&:\, u_{5,1^b}= u_{5,2^b} v_{5,2^b},\quad v_{5,1^b}= v_{5,2^b},\label{phi52b}\\
\tau_{5,2^b}&:\, u_{5,2^b}\to \frac{1}{u_{5,2^b}},\label{tau52b}\\
\phi_{5,3^b}&:\, u_{5,2^b}= -\frac{N}{t}+u_{5,3^b} v_{5,3^b},\quad v_{5,2^b}= v_{5,3^b},\label{phi53b}\\
\phi_{5,4^b}&:\, u_{5,3^b}= \frac{-1-2N+n-t+\alpha}{t}+u_{5,4^b} v_{5,4^b},\quad v_{5,3^b}= v_{5,4^b},\label{phi54b}
\end{align}
where we have denoted charts as $(u_{5,i^b}, v_{5,i^b})$ and transformations as $\phi_{5,i^b}$ for $i=1,2,3,4$ to make a distinction with the first choice of transformations. Applying this sequence of transformations we obtain the following system which is regular on the exceptional divisor $v_{5,4^b}=0$
\begin{equation}
\label{systuv424original}
\begin{gathered}
u_{5,4^b}'=-\frac{R_1^{5,4^b}(u_{5,4^b},v_{5,4^b},t)}{N t \left(v_{5,4^b} (-\alpha -n+2 N+1)+N-t u_{5,4^b} v_{5,4^b}^2\right)},\\
v_{5,4^b}'=\frac{R_2^{5,4^b}(u_{5,4^b},v_{5,4^b},t)}{N t \left(v_{5,4^b} (-\alpha -n+2 N+1)+N-t u_{5,4^b} v_{5,4^b}^2\right)},
\end{gathered}
\end{equation}
where \\

$R_1^{5,4^b}(u_{5,4^b},v_{5,4^b},t)=u_{5,4^b} (N v_{5,4^b} (-\alpha +2 n t-n+2 N+2 t+1)+N^2 (t+1)+(N+1) t v_{5,4^b}^2 (-\alpha +N+1))+N
   (N-n) (-\alpha -n+N)+N (t-1) t u_{5,4^b}^2 v_{5,4^b}^2,$\\

$R_2^{5,4^b}(u_{5,4^b},v_{5,4^b},t)=N
   v_{5,4^b}^2 (\alpha ^2-2 \alpha +n^2+2 \alpha  n+N (-4 \alpha -4 n+2 t+4)-2 n+4 N^2-2 N t u_{5,4^b}-\alpha  t+2
   t+1)+N^2 v_{5,4^b} (-2 \alpha -2 n+4 N+t+2)+t v_{5,4^b}^3 ((N+1) (-\alpha +N+1)-2 N u_{5,4^b} (-\alpha -n+2
   N+1))+N^3+N t^2 u_{5,4^b}^2 v_{5,4^b}^4$.\\

To obtain the regularising transformation from the original system~\eqref{original system2} in $(Q,P)$ coordinates to system~\eqref{systuv424original} it just remains to define the composition of the previous sequence of transformations
$$\Phi_{5,4^b}=\phi_{5,4^b}\circ\phi_{5,3^b}\circ\tau_{5,2^b}\circ\phi_{5,2^b}\circ\phi_{5,1^b},$$
which leads to
\begin{equation}
\label{Phi54}
\Phi_{5,4^b}:\ Q=\frac{t v_{5,4^b}^2}{v_{5,4^b} \left(\alpha +n-2 N+t u_{5,4^b} v_{5,4^b}-t-1\right)-N},\quad P=v_{5,4^b},
\end{equation}
and whose inverse is given by
$$ \Phi_{5,4^b}^{-1}:\, u_{5,4^b}=\frac{(-\alpha -n+2 N+t+1)P+N}{t P^2}+\frac{1}{Q},\quad v_{5,4^b}=P.$$

Although system~\eqref{systuv424original} has a more complicated structure than system~\eqref{systuv54op1}, this is also relevant since it will allow us to establish Theorem~\ref{Theorem P5}, which yields the decomposition into {\color{Brown}simpler} birational transformations of the original transformation $\hat{\phi}_{1,1}$ given by~\eqref{phi11}.

We observe that systems~\eqref{systUV11original},~\eqref{systUV21original},~\eqref{systUV31original},~\eqref{systuv54op1} in the charts $(U_{i,1},V_{i,1}),$ with $i=1,2,3$, and $(u_{5,4},v_{5,4}),$ have similar structure. Furthermore, in Section~\ref{conection PV} we will explain how these systems can be reduced to the $P_V$ equation for certain parameters. However, systems~\eqref{systUV41original} and~\eqref{systuv424original} have more complicated right--hand sides. In the following section we will iterate the regularisation procedure to find systems with simpler right--hand sides which could be easily reduced to the fifth Painlev\'{e} equation. Besides, we will show that iterated regularisation will also allow us to obtain polynomial right--hand side systems.

\section{Iterated regularisation}
\label{iterated regularisation}
Now our main objective is to identify a sequence of birational transformations that simplifies the structure of the final chart systems~\eqref{systUV41original} and~\eqref{systuv424original}. To achieve this, we can apply the same regularisation procedure employed in Section~\ref{regularisation}, which means locating the points of indeterminacy in systems~\eqref{systUV41original} and~\eqref{systuv424original}, and subsequently blowing them up until we obtain a regular system on the corresponding exceptional divisors. Among all the points of indeterminacy of these systems, we will select some cascades of points which give rise to the final chart systems with simpler structure. We iterate the procedure until we obtain systems which have analogous structure to systems~\eqref{systUV11original},~\eqref{systUV21original},~\eqref{systUV31original},~\eqref{systuv54op1}.
In Figure~\ref{diagramaIt2}, we graphically illustrate  the steps of the iterated regularisation procedure.

%%%%%%%%%%%%%%%%%%%%%%%%%%%%%%%%%%%%%%%%%%%%%%%%%%%%%%%%%%%%%%%%%%%%%%%%%%%%%%%%%%%%%   Diagram of iterated regularisation algorithm.
\begin{figure}[h!]
\begin{center}
{\footnotesize
\begin{tikzpicture}[>=latex]
\label{flowchartit2}
  %
  % Styles for states, and state edges
  %

  \tikzstyle{state} = [draw, very thick, fill=white, rectangle, minimum height=2em,minimum width=2em, node distance=5em, font={\sffamily}]
    \tikzstyle{state2} = [draw, very thick, fill=white, rectangle, minimum height=2em,text width=22em, minimum width=2em, node distance=5em, text centered, font={\sffamily}]
   % \tikzstyle{state} = [draw, very thick, fill=white, rectangle, minimum height=2em, maximum width=8em, node distance=5em, font={\sffamily}]
  \tikzstyle{stateEdgePortion} = [black,very thick];
  \tikzstyle{stateEdge} = [stateEdgePortion,->];  % sin flechas
%  \tikzstyle{stateEdge} = [stateEdgePortion,->];  % con flechas
  \tikzstyle{edgeLabel} = [pos=0.5, text centered,very thick, font={\sffamily\small}];

\node[state, name=d1] {\textbf{Initial system on $\mathbb{P}^1\times \mathbb{P}^1$}};

\node[state, name=d2,below of=d1] {Points of indeterminacy};
\draw ($(d1.south) + (0em,0)$)
      edge[stateEdge] node[edgeLabel]{}
      ($(d2.north)$);

\node[state, name=d3,below of=d2] {Cascades of points};
\draw ($(d2.south) + (0em,0)$)
      edge[stateEdge] node[edgeLabel]{}
      ($(d3.north)$);

\node[state, name=d4,below of=d3] {Final system for each cascade};
\draw ($(d3.south) + (0em,0)$)
      edge[stateEdge] node[edgeLabel]{}
      ($(d4.north)$);

\node[state2, name=d5,below of=d4] {Each system is treated as the initial system for the  next iteration and the algorithm starts again. The procedure stops when the system looks simpler, e.g. polynomial};
\draw ($(d4.south) + (0em,0)$)
      edge[stateEdge] node[edgeLabel]{}
      ($(d5.north)$);
  %\draw ($(d5.east) $)
    %  edge[stateEdge, bend right=75] node[edgeLabel]{}
     % ($(d2.east)$);

\end{tikzpicture}
}
\caption{Diagram of the iterated regularisation algorithm.}
\label{diagramaIt2}
\end{center}
\end{figure}
%%%%%%%%%%%%%%%%%%%%%%%%%%%%%%%%%%%%%%%%%%%%%%%%%%%%%%%%%%%%%%%%%%%%%%%%%%%%%%%%%%%%%%%%

As we have commented previously, Figure~\ref{diagrama} shows a detailed scheme of the general regularisation process applied to the initial system~(\ref{original system2}). We begin with the iterated regularisation for the final chart system~\eqref{systUV41original} coming from the point $P_4=\left(q=\frac{n}{N},\,P=0\right)$. For that purpose we apply the following birational transformation
\begin{align}
\label{phi42}
\phi_{4,2}&:\,\frac{1}{U_{4,1}}=  u_{4,2}v_{4,2},\quad {V_{4,1}}= v_{4,2}.
\end{align}
The resulting final chart system which is regular on the exceptional divisor $v_{4,2}=0$ is given by
\begin{equation}
\label{systuv42}
\begin{gathered}
u_{4,2}'=\frac{R_1^{4,2}(u_{4,2},v_{4,2},t)}{N t \left(n+N u_{4,2}+N v_{4,2}\right)},\\
v_{4,2}'=\frac{R_2^{4,2}(u_{4,2},v_{4,2},t)}{N t \left(n+N u_{4,2}+N v_{4,2}\right)},
\end{gathered}
\end{equation}
where \\

$R_1^{4,2}(u_{4,2},v_{4,2},t)=N u_{4,2} (-\alpha +2 n t+N^2-\alpha  N+2 N t v_{4,2}+2 N+1)+t(n+N v_{4,2}){}^2+N^3
   u_{4,2}^3+N^2 u_{4,2}^2 (-\alpha +2 N+t+2),$\\

$R_2^{4,2}(u_{4,2},v_{4,2},t)=-N v_{4,2} (\alpha -2 n^2-2
   \alpha  n+N (\alpha +4 n-2)+2 n N u_{4,2}+4 n+N^2 u_{4,2}^2-N^2-1)+N^2 v_{4,2}^2 (\alpha +n-2 N u_{4,2}-2 N-2)+n
   (-n+N+1) (-\alpha -n+N+1).$\\

\noindent At this point, we show that by applying three distinct birational transformations to system~\eqref{systuv42},  we obtain three different final charts systems which coincide with systems~\eqref{systUV11original}, \eqref{systUV21original}, \eqref{systUV31original}, respectively. We use the notation $(u_{4,3^a}, v_{4,3^a}), (u_{4,3^b}, v_{4,3^b}), (u_{4,3^c}, v_{4,3^c}) $ for the three alternative charts introduced by each transformation.

We start performing the following birational transformations to system~\eqref{systuv42}
\begin{align}
\label{phi431}
{\phi}_{4,3^a}&:\,u_{4,2}=  u_{4,3^a} v_{4,3^a},\quad u_{4,2}= -\frac{n}{N}+v_{4,3^a},
\end{align}
which leads to the following differential system
\begin{equation}
\label{systuv431}
\begin{gathered}
u_{4,3^a}'=\frac{u_{4,3^a} (-\alpha -n+2 N+t+2)-\left(u_{4,3^a}^2 \left(n-2 N v_{4,3^a}\right)\right)+t}{t},\\
v_{4,3^a}'=\frac{N V_{1,1} \left(u_{4,3^a} \left(u_{4,3^a}+2\right) \left(n-N v_{4,3 ^a}\right)+\alpha +n-2 N-2\right)+(N+1) (-\alpha +N+1)}{N t \left(u_{4,3^a}+1\right)}.
\end{gathered}
\end{equation}

This system is regular on the exceptional divisor $v_{4,3^a}=0$. As mentioned before, we notice that system~\eqref{systuv431} coincides with system~\eqref{systUV11original} in $(U_{1,1},V_{1,1})$ chart.
Alternatively to transformation ${\phi}_{4,3^a}$, we apply the following birational transformation to system~\eqref{systuv42}
\begin{align}
\label{phi432}
{\phi}_{4,3^b}&:\,u_{4,2}= -\frac{1+N}{N}+ u_{4,3^b}v_{4,3^b},\quad v_{4,2}= \frac{1+N-n}{N}+v_{4,3^b},
\end{align}
and obtain the regular system on the exceptional divisor $v_{4,3^b}=0$
\begin{equation}
\label{systuv432}
\begin{gathered}
u_{4,3^b}'=\frac{u_{4,3^b}^2 \left(-n+2 N v_{4,3^b}+N+1\right)-u_{4,3^b} (\alpha +n-t)+t}{t},\\
v_{4,3^b}'=\frac{N v_{4,3^b} \left(-u_{4,3^b} \left(u_{4,3^b}+2\right) \left(N v_{4,3 ^b}-n+N+1\right)+\alpha +n\right)+\alpha  (N+1)}{N t \left(u_{4,3^b}+1\right)}.
\end{gathered}
\end{equation}
We notice that system~\eqref{systuv432} coincides with system~\eqref{systUV21original} in $(U_{2,1},V_{2,1})$ chart.
As another alternative to the previous transformations, we now perform the following birational transformation to system~\eqref{systuv42}
\begin{align}
\label{phi433}
{\phi}_{4,3^c}&:\,u_{4,2}= \frac{-1-N+\alpha}{N}+ u_{4,3^c}v_{4,3^c},\quad u_{4,2}= \frac{1+N-n-\alpha}{N}+v_{4,3^c},
\end{align}
and obtain the following system which is regular on the exceptional divisor $v_{4,3^c}=0$
\begin{equation}
\label{systuv433}
\begin{gathered}
u_{4,3^c}'=\frac{u_{4,3^c}^2 \left(-\alpha -n+2 N v_{4,3^c}+N+1\right)+u_{4,3^c} (\alpha -n+t)+t}{t},\\
v_{4,3^c}'=\frac{\alpha  (\alpha -N-1)-N v_{4,3^c} \left(\alpha +u_{4,3^c} \left(U_{4,3 ^c}+2\right) \left(-\alpha -n+N v_{4,3^c}+N+1\right)-n\right)}{N t \left(u_{4,3 ^c}+1\right)}.
\end{gathered}
\end{equation}
Once again, we observe that system~\eqref{systuv433} coincides with system~\eqref{systUV31original} in $(U_{3,1},V_{3,1})$ chart. The previous observations can be summarised in the following statement.\\

\begin{theorem}
\label{Theorem P4}
System~\eqref{systuv42} in $(u_{4,2}, v_{4,2})$ chart can be transformed to systems~\eqref{systUV11original},~\eqref{systUV21original},~\eqref{systUV31original} in $(U_{1,1}V_{1,1})$, $(U_{2,1}V_{2,1})$, $(U_{3,1}V_{3,1})$ charts  by the following birational transformations
\begin{align*}
&\hat{{\varphi}}_{1,1}:\, u_{4,2}=  U_{1,1}V_{1,1},\quad v_{4,2}= -\frac{n}{N}+V_{1,1},\\
&\hat{{\varphi}}_{2,1}: u_{4,2}= -\frac{1+N}{N}+ U_{2,1}V_{2,1},\quad v_{4,2}= \frac{1+N-n}{N} +V_{2,1},\\
&\hat{{\varphi}}_{3,1}: u_{4,2}= \frac{-1-N+\alpha}{N} +U_{3,1}V_{3,1},\quad v_{4,2}= \frac{1+N-n-\alpha}{N}+V_{3,1},
\end{align*}
respectively. Furthermore, iterated regularisation leads to the following decompositions:
\begin{align*}
\hat{\phi}_{1,1}=\phi_{q,P}\circ\hat{{\varphi}}_{1,1}\circ\phi_{4,3^a}\circ\phi_{4,2}\circ\hat{\phi}_{4,1},\\
\hat{\phi}_{2,1}=\phi_{q,P}\circ\hat{{\varphi}}_{2,1}\circ\phi_{4,3^b}\circ\phi_{4,2}\circ\hat{\phi}_{4,1},\\
\hat{\phi}_{3,1}=\phi_{q,P}\circ\hat{{\varphi}}_{3,1}\circ\phi_{4,3^c}\circ\phi_{4,2}\circ\hat{\phi}_{4,1},
\end{align*}
where transformations $\hat{\phi}_{i,1},\, i=1,2,3,4$ are defined by~\eqref{phi11},~\eqref{phi21},~\eqref{phi31},~\eqref{phi41}, $\phi_{4,2}$ is given by~\eqref{phi42}, $\phi_{4,3^a}, \phi_{4,3^b}, \phi_{4,3^c}$ are given by~\eqref{phi431},~\eqref{phi432},~\eqref{phi433}, and $\phi_{q,P}$ is defined by
\begin{equation*}
\phi_{q,P}:\, q=q,\quad P=\frac{1}{p}.
\label{phiqP}
\end{equation*}
\end{theorem}

Similarly, we proceed with iterated regularisation for final chart system~\eqref{systuv424original} coming from the point $P_5=(Q=0,\,P=0)$. We consider system~\eqref{systuv424original} and apply a second iteration of the regularisation procedure. We take the birational transformation
\begin{equation}
\label{phi55}
\phi_{5,5^b}:\,u_{5,4^b}= u_{5,5^b} v_{5,5^b},\quad \frac{1}{v_{5,4^b}}= v_{5,5^b},
\end{equation}
which yields the final chart system which is regular on the exceptional divisor $v_{5,5^b}=0$
\begin{equation}
\label{systuv55}
\begin{gathered}
u_{5,5^b}'=\frac{R_1^{5,5^b}(u_{5,5^b},v_{5,5^b},t)}{t \left(\alpha +n-N v_{5,5^b}-2 N+t u_{5,5^b}-1\right)},\\
v_{5,5^b}'=\frac{R_2^{5,5^b}(u_{5,5^b},v_{5,5^b},t)}{N t \left(\alpha +n-N v_{5,5^b}-2 N+t u_{5,5^b}-1\right)},
\end{gathered}
\end{equation}
where \\

$R_1^{5,5^b}(u_{5,5^b},v_{5,5^b},t)=-u_{5,5^b} (\alpha ^2-\alpha +n^2+2 \alpha  n+N (-4 \alpha -4 n+2 t+2)+N v_{5,5^b} (-2 \alpha -2
   n+4 N+1)-2 n t-n+N^2 v_{5,5^b}^2+4 N^2-\alpha  t)+(N-n) (-\alpha -n+N)+t u_{5,5^b}^2 (-2 \alpha -2 n+2 N
   v_{5,5^b}+4 N+t+1)-t^2 u_{5,5^b}^3,$\\

$R_2^{5,5^b}(u_{5,5^b},v_{5,5^b},t)=N v_{5,5^b}
   (\alpha ^2-2 \alpha +n^2+2 \alpha  n+N (-4 \alpha -4 n+2 t+4)-2 t u_{5,5^b} (-\alpha -n+2 N+1)-2 n+4 N^2+t^2
   u_{5,5^b}^2-\alpha  t+2 t+1)+N^2 v_{5,5^b}^2 (-2 \alpha -2 n+4 N-2 t u_{5,5^b}+t+2)+N^3 v_{5,5^b}^3+(N+1) t
   (-\alpha +N+1).$\\

Next, we apply a third iteration of the regularisation procedure. We consider system~\eqref{systuv55} and apply the birational transformation
\begin{equation}
\label{phi56}
\phi_{5,6^b}:\,{u_{5,5^b}}= u_{5,6^b} v_{5,6^b},\quad \frac{1}{v_{5,5^b}}=  v_{5,6^b},
\end{equation}
obtaining the following system regular on the exceptional divisor $v_{5,6^b}=0$
\begin{equation}
\label{systuv56}
\begin{gathered}
u_{5,6^b}'=-\frac{R_1^{5,6^b}(u_{5,6^b},v_{5,6^b},t)}{N t \left(N-v_{5,6^b} \left(\alpha +n-2 N+t u_{5,6^b}
   v_{5,6^b}-1\right)\right)},\\
v_{5,6^b}'=\frac{R_2^{5,6^b}(u_{5,6^b},v_{5,6^b},t)}{N t
   \left(N-v_{5,6^b} \left(\alpha +n-2 N+t u_{5,6^b} v_{5,6^b}-1\right)\right)},
\end{gathered}
\end{equation}
where \\

$R_1^{5,6^b}(u_{5,6^b},v_{5,6^b},t)=u_{5,6^b} (N v_{5,6^b} (-\alpha +2 n t-n+2 N+2 t+1)+N^2 (t+1)+t v_{5,6^b}^2 ((N+1) (-\alpha
   +N+1)+N (t-1) u_{5,6^b}))+N (N-n) (-\alpha -n+N),$\\

$R_2^{5,6^b}(u_{5,6^b},v_{5,6^b},t)=v_{5,6^b} (v_{5,6^b} (N (\alpha ^2+n^2+2 (\alpha -1) n+2 N (-2
   \alpha -2 n+t+2)+4 N^2-\alpha  (t+2)+2 t+1)+t (2 N u_{5,6^b}(v_{5,6^b} (\alpha +n-2 N-1)-N)+N t
   u_{5,6^b}^2 v_{5,6^b}^2+(N+1) v_{5,6^b} (-\alpha +N+1)))+N^2 (-2 \alpha -2 n+4 N+t+2))+N^3.$\\

Now, we apply a fourth iteration of the regularisation procedure, which involves performing the following birational sequence of transformations to system~\eqref{systuv56}
\begin{align}
\phi_{5,7^b}&:\,\frac{1}{u_{5,6^b}}= u_{5,7^b} v_{5,7^b},\quad v_{5,6^b}=v_{5,7^b},\label{phi57b}\\
\phi_{5,8^b}&:\,{u_{5,7^b}}= u_{5,8^b} v_{5,8^b},\quad v_{5,7^b}=v_{5,8^b},\label{phi58b}\\
\tau_{5,8^b}&:\,u_{5,8^b}\to \frac{1}{u_{5,8^b}},\label{tau58b}\\
\phi_{5,9^b}&:\,u_{5,8^b}= \frac{N}{t}+u_{5,9^b} v_{5,9^b},\quad v_{5,8^b}= v_{5,9^b},\label{phi59b}\\
\phi_{5,10^b}&:\,u_{5,9^b}= \frac{1+2N-n+t+\alpha}{t}+u_{5,10^b} v_{5,10^b},\quad v_{5,9^b}= v_{5,10^b},\label{phi510b}
\end{align}
to obtain the following system which is regular on the exceptional divisor $v_{5,10^b}=0$
\begin{equation}
\label{systuv510}
\begin{gathered}
u_{5,10^b}'=\frac{R_1^{5,10^b}(u_{5,10^b},v_{5,10^b},t)}{N t v_{5,10^b} \left(u_{5,10^b}
   v_{5,10^b}+1\right)},\\
v_{5,10^b}'=-\frac{R_2^{5,10^b}(u_{5,10^b},v_{5,10^b},t)}{N t \left(u_{5,10^b} v_{5,10^b}+1\right)},
\end{gathered}
\end{equation}
where \\

$R_1^{5,10^b}(u_{5,10^b},v_{5,10^b},t)=u_{5,10^b} (2 n N v_{5,10^b}-N^2+(N+1) v_{5,10^b}^2 (-\alpha +N+1))+N u_{5,10^b}^2 v_{5,10^b}(v_{5,10^b} (\alpha +n-2 N-2)-2 N)+n N,$\\

$R_2^{5,10^b}(u_{5,10^b},v_{5,10^b},t)=N^2+v_{5,10^b}^2 ((N+1) (-\alpha +N+1)+2 N t u_{5,10^b})+N t u_{5,10^b}^2
   v_{5,10^b}^3+N v_{5,10^b} (-\alpha +2 N+t+2).$\\

\noindent Furthermore, the regularising transformation from system~\eqref{systuv56} to system~\eqref{systuv510} is deduced from the composition of the previous transformations, that is,
$$\Phi_{5,10^b}=\phi_{5,10^b}\circ\phi_{5,9^b}\circ\tau_{5,8^b}\circ\phi_{5,8^b}\circ\phi_{5,7^b}, $$
\begin{equation}
\label{Phi510}
\Phi_{5,10^b}:\, \frac{1}{u_{5,6^b}}= \frac{t v_{5,10^b}^2}{v_{5,10^b} \left(-\alpha -n+2 N+t u_{5,10^b} v_{5,10^b}+t+1\right)+N},\quad v_{5,6^b}= v_{5,10^b},
\end{equation}
and its inverse is given by
 $$\Phi_{5,10^b}^{-1}:\, u_{5,10^b}=u_{5,6^b}+\frac{-N+(-1-2N+n-t+\alpha)v_{5,6^b}}{tv_{5,6^b}^2} ,\quad v_{5,10^b}= v_{5,6^b}.$$
Finally, we perform a fifth iteration of the regularisation procedure. We observe that applying another blow--up to system~\eqref{systuv510} we return to system~\eqref{systUV11original} in the original chart $(U_{1,1},V_{1,1})$. This result is summarised in the following statement.\\

\begin{theorem}
\label{Theorem P5}
System~\eqref{systuv510} in $(u_{5,{10}^b}, v_{5,{10}^b})$ chart is transformed to system~\eqref{systUV11original} in $(U_{1,1},V_{1,1})$ chart  by the transformation
\begin{equation}
{\hat{\psi}}_{1,1}:\, u_{5,10^b}=  V_{1,1},\quad v_{5,10^b}= U_{1,1}V_{1,1},\label{psi11}
\end{equation}
which gives rise to the following decomposition:
$$\hat{\phi}_{1,1}=\phi_{Q,P}\circ{\hat{\psi}}_{1,1}\circ\Phi_{5,10^b}\circ{\phi}_{5,6^b}\circ\phi_{5,5^b}\circ\Phi_{5,4^b},$$
where $\hat{\phi}_{1,1}$ is given by~\eqref{phi11}, $\phi_{5,i^b},\ i=5,6$, $\Phi_{5,j^b},\ j=4,10$, are given by~\eqref{Phi54},~\eqref{phi55},~\eqref{phi56},~\eqref{Phi510} and $\phi_{Q,P}$
is defined by
\begin{equation*}
\label{phiQP}
\phi_{Q,P}:\, Q=\frac{1}{q},\quad P=\frac{1}{p}.
\end{equation*}
\end{theorem}

Now, our main objective is to investigate the connection between systems~\eqref{systUV11original}, \eqref{systUV21original}, \eqref{systUV31original}, \eqref{systuv54op1} and the fifth Painlev\'{e} equation. We use two approaches to find this connection. First we reduce systems~\eqref{systUV11original}, \eqref{systUV21original}, \eqref{systUV31original}, \eqref{systuv54op1} to second--order differential equations which, after applying a M\"{o}bius transformation, become $P_V$ for certain values of the parameters. The second alternative consists in performing an additional blow--up to systems~\eqref{systUV11original}, \eqref{systUV21original}, \eqref{systUV31original}, \eqref{systuv54op1}. This leads to polynomial right--hand side systems, which can be also reduced to second--order differential equations that directly coincide with the $P_V$ equation with the same parameters as previously had been obtained for each system, without requiring any further transformation.

\subsection{Connection to the fifth Painlev\'{e} equation}
\label{conection PV}
This section is devoted to study the connection between systems~\eqref{systUV11original},~\eqref{systUV21original},~\eqref{systUV31original}, \eqref{systuv54op1} and the fifth Painlev\'{e} equation given by~\eqref{PV}. In Figure~\ref{diagrama} we give a graphical representation of the regularisation procedure and show how we establish connection with the fifth Painlev\'{e} equation.

To find this connection we begin observing that these systems can be reduced to second--order differential equations with the $U_{i,1}$ (for $i=1,2,3$) and $v_{5,4}$ variables. For systems \eqref{systUV11original}, \eqref{systUV21original}, \eqref{systUV31original} we solve for $V_{i,1}$ in the first equation of each system and substitute this value into the second equation. Analogously, for system~\eqref{systuv54op1} we solve for $u_{5,4}$ in the second equation and substitute this value into the first equation.  For instance, solving for $V_{1,1}$ in the first equation of system~\eqref{systUV11original} we find
\begin{equation}
\label{V11}
V_{1,1}= \frac{t \left(U_{1,1}'-1\right)+U_{1,1} (\alpha +n-2 N-t-2)+n U_{1,1}^2}{2 N U_{1,1}^2},
\end{equation}
and substituting this value into the second equation of system~\eqref{systUV11original},  the following second--order differential equation is obtained

\begin{align}
U_{1,1}''=& \frac{\left(\frac{1}{U_{1,1}}+\frac{3}{2}\right) \left(U_{1,1}'\right)^2}{U_{1,1}+1}-\frac{U_{1,1}'}{t}+\frac{U_{1,1}^2
   \left(-\alpha +n U_{1,1}+n\right) \left(\alpha +n U_{1,1}+n\right)}{2 t^2 \left(U_{1,1}+1\right)}\notag\\
   &+\frac{U_{1,1}^2 (\alpha +n-2 N-1)}{t \left(U_{1,1}+1\right)}
   +\frac{2 \alpha +4 U_{1,1} (\alpha +n-2 N-1)+2 n-4 N-5 t-2}{2 t \left(U_{1,1}+1\right)}\notag\\
   &-\frac{\left(U_{1,1}+4\right) U_{1,1}}{2 \left(U_{1,1}+1\right)}-\frac{1}{\left(U_{1,1}+1\right) U_{1,1}}.\label{equationU11}
\end{align}

Proceeding in an analogous way for the remaining systems we deduce the following second--order differential equations:

\begin{align}
U_{2,1}''=& \frac{\left(\frac{1}{U_{2,1}}+\frac{3}{2}\right) \left(U_{2,1}'\right)^2}{U_{2,1}+1}-\frac{ U_{2,1}'}{t}\notag\\
&+\frac{U_{2,1}^2 \left((n-\alpha ) (\alpha +n-2 N-2)+U_{2,1} \left(U_{2,1}+2\right) (-n+N+1)^2\right)}{2 t^2
   \left(U_{2,1}+1\right)}+\frac{U_{2,1}^2 (\alpha +n+1)}{t \left(U_{2,1}+1\right)}\notag\\
   &+\frac{2 \alpha +4 U_{2,1} (\alpha +n+1)+2 n-5 t+2}{2 t \left(U_{2,1}+1\right)}-\frac{\left(U_{2,1}+4\right) U_{2,1}}{2 \left(U_{2,1}+1\right)}-\frac{1}{\left(U_{2,1}+1\right) U_{2,1}}\label{equationU21},
\end{align}

\begin{align}
U_{3,1}''=& \frac{\left(\frac{1}{U_{3,1}}+\frac{3}{2}\right) \left(U_{3,1}'\right)^2}{U_{3,1}+1}-\frac{ U_{3,1}'}{t }\notag\\
&+\frac{U_{3,1}^2 \left((\alpha +n) (\alpha +n-2 N-2)+U_{3,1} (U_{3,1}+2)(\alpha +n-N-1)^2 \right)}{2 t^2
   \left(U_{3,1}+1\right)}+\frac{U_{3,1}^2 (-\alpha +n+1)}{t \left(U_{3,1}+1\right)}\notag\\
&+\frac{-2 \alpha +4 U_{3,1} (-\alpha +n+1)+2 n-5 t+2}{2 t \left(U_{3,1}+1\right)}-\frac{(U_{3,1}+4)U_{3,1}}{2 \left(U_{3,1}+1\right)}-\frac{1}{\left(U_{3,1}+1\right) U_{3,1}}\label{equationU31},
\end{align}

\begin{align}
v_{5,4}''=& \frac{\left(\frac{1}{v_{5,4}}+\frac{3}{2}\right) \left(v_{5,4}'\right)^2}{v_{5,4}+1}-\frac{v_{5,4}'}{t}+\frac{v_{5,4}^2
   \left(-\alpha +n v_{5,4}+n\right) \left(\alpha +n v_{5,4}+n\right)}{2 t^2 \left(v_{5,4}+1\right)}\notag\\
   &+\frac{v_{5,4}^2 (\alpha +n-2 N-1)}{t \left(v_{5,4}+1\right)}+\frac{2 \alpha +4 v_{5,4} (\alpha +n-2 N-1)+2 n-4 N-5 t-2}{2 t \left(v_{5,4}+1\right)}\notag\\
   &-\frac{\left(v_{5,4}+4\right) v_{5,4}}{2 \left(v_{5,4}+1\right)}-\frac{1}{\left(v_{5,4}+1\right) v_{5,4}}.\label{equationv54}
\end{align}

We observe that equation~\eqref{equationv54} coincides with equation~\eqref{equationU11} renaming the variables. Finally, it is remarkable that these equations exhibit, choosing the suitable parameters, a similar structure to equation $P_V$ given in~\eqref{PV}. In particular, applying the M\"{o}bius transformation

\begin{equation}
\label{transf moebius}
U_{1,1}(t)=U_{2,1}(t)=U_{3,1}(t)=-1+y(t),
\end{equation}
in equation~\eqref{equationU11}, we obtain $P_V$ with parameters:
\begin{equation}
A_5=\frac{n^2}{2},\,B_5=-\frac{\alpha ^2}{2},\,C_5=\alpha+n -2 N-1,\,D_5=-\frac{1}{2}.\label{param1}
\end{equation}
Similarly, for equation~\eqref{equationU21} we deduce $P_V$ with
\begin{equation}
\label{param2}
A_5=\frac{1}{2} (N-n+1)^2,\,B_5=-\frac{1}{2} (-\alpha +N+1)^2,\,C_5=\alpha +n+1,\,D_5=-\frac{1}{2},
\end{equation}
and for equation~\eqref{equationU31} we have
\begin{equation}
A_5=\frac{1}{2} (N-n+1-\alpha)^2,\,B_5=-\frac{1}{2} (1+N)^2,\,C_5=-\alpha+n+1,\,D_5=-\frac{1}{2}.\label{param3}
\end{equation}

Now, we examine how we can find relations between $P_V$ equation~\eqref{PV} with different parameters~\eqref{param1},~\eqref{param2} and~\eqref{param3}. For that purpose, since $D_5\neq0$, the B\"{a}cklund transformation can be written as (see~\cite[Th. 39.2]{Gromak})
\begin{equation}
\label{backlund}
T_{\varepsilon_1,\varepsilon_2,\varepsilon_3}: y=y(t,A_5,B_5,C_5,D_5)\to y_1=y_1(t,A_5^{(1)},B_5^{(1)},C_5^{(1)},D_5^{(1)}),
\end{equation}
where
\newpage
$$y_1=1-\frac{2\varepsilon_3kty}{ty'-\varepsilon_1cy^2+(\varepsilon_1 c-\varepsilon_2a+\varepsilon_3kt)y+\varepsilon_2a},$$
with $c^2=2A_5$, $a^2=-2B_5$, $k^2=-2D_5$, $\varepsilon_1^2=\varepsilon_2^2=\varepsilon_3^2=1$ and
\begin{align*}
A_5^{(1)}&=-\frac{1}{16D_5}(C_5+\varepsilon_3k(1-\varepsilon_2a-\varepsilon_1c))^2,\\
B_5^{(1)}&=\frac{1}{16D_5}(C_5-\varepsilon_3k(1-\varepsilon_2a-\varepsilon_1c))^2,\\
C_5^{(1)}&=\varepsilon_3k(\varepsilon_2a-\varepsilon_1c),\\
D_5^{(1)}&=D_5.
\end{align*}

In the following theorem we show that compositions of B\"{a}cklund transformations~\eqref{backlund}, which may not be unique, allow us to establish relations between $P_V$ equation~\eqref{PV} with parameters~\eqref{param1},~\eqref{param2} and~\eqref{param3}. We first note that in~\eqref{param1},~\eqref{param2} and~\eqref{param3} we have $D_5=-\frac{1}{2}$, thus, $k=1$. To obtain the following B\"{a}cklund  transformations we fix positive branches in parameters $c$ and $a$, that is, $c=\sqrt{2A_5}$ and $a=\sqrt{-2B_5}.$\\

\begin{theorem}The fifth Painlev\'{e} equation~\eqref{PV} with parameters~\eqref{param1} is related to equation~\eqref{PV} with parameters~\eqref{param2} through the following composition of B\"{a}cklund  transformations
$$T_{-1,-1,-1}\circ T_{1,1,-1}\circ T_{1,1,1}\circ T_{1,-1,1},$$
given by~\eqref{backlund}, which gives rise to the transformation
$$y\to y-\frac{2 (N+1) (y-1)^2 y}{y (-\alpha +n-2 N+t-2)+y^2 (-n+2 N+2)-t y'+\alpha}.$$
Similarly, the composition
$$T_{1,1,-1}\circ T_{1,1,1}\circ T_{1,1,-1}\circ T_{1,1,1},$$ which produces the following transformation
$$y\to y-\frac{2 (y-1)^2 y (-\alpha +N+1)}{y (3 \alpha +n-2 N+t-2)+y^2 (-2 \alpha -n+2 N+2)-t y'-\alpha },$$
converts the parameters~\eqref{param1} to the parameters~\eqref{param3}.
Finally, the composition
\newpage
 $$T_{1,1,1}\circ T_{1,-1,-1},$$ which leads to the transformation
 $$y\to y+\frac{2 \alpha  y (y-1)^2}{y (3 \alpha +n-2 N+t-2)+y^2 (-2 \alpha -n+N+1)+N-t y'+1-\alpha},$$
converts parameters~\eqref{param2} to~\eqref{param3}.
\end{theorem}

\begin{remark} Transformations $U_{1,1}(t)=U_{2,1}(t)=U_{3,1}(t)=-1+\frac{1}{y(t)}$ in
equations~\eqref{equationU11}, \eqref{equationU21}, \eqref{equationU31} also result in $P_V$ for some alternative parameters. For instance, by using the previous transformation in equation~\eqref{equationU11}, $P_V$ is obtained for parameters $A_5=\frac{\alpha^2}{2}, B_5=-\frac{n^2}{2},C_5=1+2N-n-\alpha,D_5=-\frac{1}{2}.$
\end{remark}

\begin{remark} In case $\alpha=0$, instead of analysing system~\eqref{systUV21original} we work with system~\eqref{systUV220} which is regular on the exceptional divisor $\tilde{V}_{2,2}=0$. From this system we deduce the following second--order differential equation for $\tilde{V}_{2,2}$

\begin{align*}
\tilde{V}_{2,2}''=& \frac{\left(1-\frac{1}{2 \tilde{V}_{2,2}}\right) \left(\tilde{V}_{2,2}'\right)^2}{\tilde{V}_{2,2}-1}-\frac{\tilde{V}_{2,2}'}{t}+\frac{n \tilde{V}_{2,2}^2 (-n+2 N+2)-2 (N+1)^2 \tilde{V}_{2,2}+(N+1)^2}{2 t^2 \left(\tilde{V}_{2,2}-1\right) \tilde{V}_{2,2}}\notag\\
&+\frac{2 (n+1) \tilde{V}_{2,2}^2}{t \left(\tilde{V}_{2,2}-1\right)}
+\frac{\tilde{V}_{2,2} \left(\tilde{V}_{2,2} \left(\tilde{V}_{2,2} \left(2 t \tilde{V}_{2,2}-2 n-5
   t-2\right)+4 t\right)-2 n-t-2\right)}{2 t \left(\tilde{V}_{2,2}-1\right)},
\end{align*}
which can be transformed to Painlev\'{e} V applying the transformation $\tilde{V}_{2,2}(t)=\frac{y(t)}{y(t)-1}$
for the parameters
$$A_5=\frac{1}{2} (N-n+1)^2,\,B_5=-\frac{1}{2} (N+1)^2,\,C_5=n+1,\,D_5=-\frac{1}{2}.$$
Notice that these parameters coincide with the previous parameters given by~\eqref{param2} when $\alpha=0.$
\end{remark}

\begin{remark} As commented previously, in~\cite{Krawtchouk Polynomials} the authors obtained some connection between system~\eqref{original system2} and $P_V$ involving some cumbersome expressions. In particular, they obtained connection with $P_V$ for the parameters
\begin{equation}
\label{parameters walter}
A_5=\frac{(\alpha-N-1)^2}{2},\quad B_5=-\frac{(n-N)^2}{2},\quad C_5=-(n+\alpha), \quad D_5=-\frac{1}{2}.
\end{equation}
Equation $P_V$ with parameters~\eqref{param1} can be connected to $P_V$ with parameters~\eqref{parameters walter} via the following B\"{a}cklund transformation $$T_{-1,1,-1}: y\to\frac{2 t y}{\alpha -y (\alpha +n+t)+n y^2+t y'}+1.$$

\end{remark}

\subsection{Iterated regularisation for obtaining polynomial right--hand side systems}
\label{polynomial reg}
Previously, we have seen that two and four more iterations of the regularisation procedure for systems~\eqref{systUV41original} and~\eqref{systuv424original}, respectively, gives rise to systems having simpler structures which can be connected to $P_V$. Moreover, applying iterated regularisation in systems~\eqref{systUV41original} and~\eqref{systuv424original} we obtained system~\eqref{systUV11original},~\eqref{systUV21original} and~\eqref{systUV31original}, as seen in Theorems~\ref{Theorem P4} and~\ref{Theorem P5}. To conclude this section we will show that a new iteration of the regularisation procedure gives rise to polynomial right--hand side systems. Specifically, we will show that systems \eqref{systUV11original}, \eqref{systUV21original}, \eqref{systUV31original}, \eqref{systuv54op1} can be reduced to polynomial right--hand side systems doing just another blow--up. These final chart systems are of special interest since they allow us to find direct connection with $P_V$. That is, the second--order differential equations for the variables $V_{i,2},\,i=1,2,3 $ and $u_{5,4}$ coincide with $P_V$ without the need of applying any additional transformation. Moreover, we will prove that these polynomials systems are also Hamiltonian.

We begin with the following birational transformation
\begin{equation}
\label{phi12}
\hat{\phi}_{1,2}:\,U_{1,1}= -1+V_{1,2},\quad V_{1,1}= \frac{1+N}{N}+ U_{1,2}V_{1,2},
\end{equation}
to system~\eqref{systUV11original} in $(U_{1,1},V_{1,1})$ chart which leads to the polynomial right--hand side system

\begin{equation}
\begin{gathered}
\label{systUV12}
{N t} U_{1,2}'=-N U_{1,2} (V_{1,2} (-2 n+4 N+4)-\alpha +n-2 N+t-2)-N^2 U_{1,2}^2 (3
   V_{1,2}^2+4 V_{1,2}-1)\\
   \hspace{-8.2cm}-(N+1) (-n+N+1),\\
 V_{1,2}'=\frac{V_{1,2} (2 N U_{1,2}-\alpha +n-2
   N+t-2)-V_{1,2}^2 (4 N U_{1,2}+n-2 N-2)+2 N U_{1,2} V_{1,2}^3+\alpha}{t}.
\end{gathered}
\end{equation}

For system~\eqref{systUV21original} in $(U_{2,1},V_{2,1})$ chart we perform the birational transformation
\begin{equation}
\label{phi22}
\hat{\phi}_{2,2}:\,U_{2,1}= -1+V_{2,2},\quad V_{2,1}= -\frac{1+N}{N}+ U_{2,2}V_{2,2},
\end{equation}
and obtain the system
\begin{equation}
\begin{gathered}
\label{systUV22}
 NtU_{2,2}'=N U_{2,2} \left(2 V_{2,2} (n+N+1)+\alpha -n-2 N-t-2\right)-N^2 U_{2,2}^2 \left(3 V_{2,2}^2-4
   V_{2,2}+1\right)\\
   -n N-n,\hspace{8.5cm}\\
 tV_{2,2}'=V_{2,2} (2 N U_{2,2}-\alpha +n+2
   N+t+2)-V_{2,2}^2 (4 N U_{2,2}+n+N+1)+2 N U_{2,2} V_{2,2}^3\\
   +\alpha -N-1,\hspace{8.5cm}
\end{gathered}
\end{equation}

which also has the polynomial right--hand side.\\

Similarly, for system~\eqref{systUV31original} in $(U_{3,1},V_{3,1})$ chart we apply the birational transformation
\begin{equation}
\label{phi32}
\hat{\phi}_{3,2}:\,U_{3,1}= -1+V_{3,2},\quad V_{3,1}= \frac{\alpha}{N}+ U_{3,2}V_{3,2},
\end{equation}

and arrive to the polynomial right--hand side system
\begin{align}
N t U_{3,2}'=&-(N U_{3,2} (2 V_{3,2} (\alpha -n+N+1)-\alpha +n-2 N+t-2)\notag\\
&+N^2 U_{3,2}^2 (3 V_{3,2}^2-4)
   V_{3,2}+1)+\alpha  (-n+N+1)),\notag\\\notag\\
 tV_{3,2}'=&V_{3,2} (2 N U_{3,2}-\alpha +n-2N+t-2)\notag\\
 &+V_{3,2}^2 (-4 N U_{3,2}+\alpha -n+N+1+2 N U_{3,2} V_{3,2}^3+N+1).\label{systUV32}
\end{align}

Analogously, for system~\eqref{systuv54op1} in $(u_{5,4},v_{5,4})$ chart we apply the birational transformation
\begin{equation}
\label{phi55}
{\phi}_{5,5}:\,u_{5,4}= -1+\frac{n}{N}+u_{5,5}v_{5,5},\quad v_{5,4}= -1+v_{5,5},
\end{equation}
leading to the polynomial right--hand side system
\begin{equation}
\begin{gathered}
\label{systuv552}
u_{5,5}'=\frac{u_{5,5} \left(v_{5,5} (4 N-2 n)-\alpha +n-2 N+t\right)-N u_{5,5}^2 \left(3 v_{5,5}^2-4
   v_{5,5}+1\right)+n-N}{t},\\
v_{5,5}'=\frac{v_{5,5} \left(2 N u_{5,5}+\alpha -n+2 N-t\right)+v_{5,5}^2 \left(-4 N u_{5,5}+n-2 N\right)+2 N
   u_{5,5} v_{5,5}^3-\alpha }{t}.
\end{gathered}
\end{equation}

Now, we see that solving for $U_{i,2},\,i=1,2,3$ and $u_{5,5}$ in the second equation of each system and substituting the obtained expressions into the first equation of the corresponding system we obtain a differential equation for $V_{i,2}$ and $v_{5,5}$ variables, respectively. For instance, solving for $U_{1,2}$ in the second equation of~\eqref{systUV12} we obtain
\begin{equation}
\label{U12}
U_{1,2}=\frac{V_{1,2} (\alpha -n+2 N-t+2)+V_{1,2}^2 (n-2 N-2)+t V_{1,2}'-\alpha }{2 N \left(V_{1,2}-1\right){}^2 V_{1,2}},
\end{equation}
and substituting this value into the first equation of the system we obtain
\begin{align}
V_{1,2}''=& \left(\frac{1}{2 V_{1,2}}+\frac{1}{V_{1,2}-1}\right) \left(V_{1,2}'\right)^2-\frac{V_{1,2}'}{t}+\frac{\left(V_{1,2}-1\right){}^2 \left(\frac{1}{2} n^2 V_{1,2}-\frac{\alpha ^2}{2 V_{1,2}}\right)}{t^2}\notag\\
&+\frac{V_{1,2} (\alpha +n-2 N-1)}{t}-\frac{V_{1,2} \left(V_{1,2}+1\right)}{2
   \left(V_{1,2}-1\right)}\label{eq12}.
\end{align}

Similarly for systems~\eqref{systUV22},~\eqref{systUV32} we obtain the following second--order differential equations

\begin{align}
V_{2,2}''=& \left(\frac{1}{2 V_{2,2}}+\frac{1}{V_{2,2}-1}\right) \left(V_{2,2}'\right)^2-\frac{V_{2,2}'}{t}\notag\\
&+\frac{\left(V_{2,2}-1\right){}^2 \left(\frac{1}{2} V_{2,2} (-n+N+1)^2-\frac{(-\alpha +N+1)^2}{2 V_{2,2}}\right)}{t^2}+\frac{V_{2,2} (\alpha +n+1)}{t}\notag\\
&-\frac{V_{2,2} \left(V_{2,2}+1\right)}{2
   \left(V_{2,2}-1\right)}\label{eq2.1},\\\notag\\
V_{3,2}''=& \left(\frac{1}{2 V_{3,2}}+\frac{1}{V_{3,2}-1}\right) \left(V_{3,2}'\right)^2-\frac{V_{3,2}'}{t}\notag\\
&+\frac{\left(V_{3,2}-1\right){}^2 \left(\frac{1}{2} V_{3,2} (-\alpha -n+N+1)^2+\frac{-N^2-2 N-1}{2 V_{3,2}}\right)}{t^2}+\frac{V_{3,2} (-\alpha +n+1)}{t}\notag\\
&-\frac{V_{3,2} \left(V_{3,2}+1\right)}{2
   \left(V_{3,2}-1\right)}.\label{eq31}
\end{align}
The equation obtained in $v_{5,5}$ variable is the same as equation~\eqref{eq12}. Furthermore, equations~\eqref{eq12},~\eqref{eq2.1},~\eqref{eq31} coincide with $P_V$ equation~\eqref{PV} for parameters~\eqref{param1},~\eqref{param2},~\eqref{param3}, respectively.\\

To conclude this section we present a diagram representing graphically the regularisation procedure and the different ways to establish a connection with the fifth Painlev\'{e} equation.

%%%%%%%%%%%%%%%%%%%%%%%%%%%%%%%%%%%%%%%%%%%%%%%%%%%%%%%%%%%%%%%%%%%%%%%%%%%%%%%%%%%%%% Flow diagram

\begin{figure}[H]
{\footnotesize
\hspace{-2.25cm}
\begin{tikzpicture}[>=latex]
\label{flowchart}
  %
  % Styles for states, and state edges
  %
  \tikzstyle{state} = [draw, very thick, fill=white, rectangle, minimum height=2em, minimum width=2em, node distance=4em, font={\sffamily}]
  \tikzstyle{statep} = [draw, very thick, fill=white, rectangle, minimum height=2em, minimum width=2em, node distance=6em, font={\sffamily}]
  \tikzstyle{statepred} = [draw=red, very thick, fill=white, rectangle, minimum height=2em, minimum width=2em, node distance=6em, font={\color{red}{\sffamily}}]
  \tikzstyle{state1} = [draw, very thick, fill=white, rectangle, minimum height=2em,text width=15em, minimum width=2em, node distance=5em, font={\sffamily}]
  \tikzstyle{state11} = [draw, very thick, fill=white, rectangle, text width=1.5cm, minimum width=1em, minimum height=2em, node distance=6em, text centered,font={\sffamily}]
  \tikzstyle{state11b} = [draw, very thick, fill=white, rectangle, text width=1.2cm, minimum width=1em,  minimum height=2em, node distance=6em, text centered,font={\sffamily}]
  \tikzstyle{state11c} = [draw, very thick, fill=white, rectangle, text width=1.5cm, minimum width=1em,  minimum height=2em, node distance=6em, text centered,font={\sffamily}]
  \tikzstyle{state11c1} = [draw, very thick, fill=white, rectangle, text width=1.5cm, minimum width=1em,  minimum height=2em, node distance=15em, text centered,font={\sffamily}]
  \tikzstyle{state11d} = [draw, very thick, fill=white, rectangle, text width=0.5cm, minimum width=1em,  minimum height=2em, node distance=6em, text centered,font={\sffamily}]
  \tikzstyle{state11e} = [draw, very thick, fill=white, rectangle, text width=3cm, minimum width=1em, minimum height=2em, node distance=12em, text centered,font={\sffamily}]
  \tikzstyle{state11f} = [draw, very thick, fill=white, rectangle, text width=2.5cm, minimum width=1em, minimum height=2em, node distance=12em, text centered,font={\sffamily}]
  \tikzstyle{stateEdgePortion} = [black,very thick];
  \tikzstyle{stateEdge} = [stateEdgePortion];  % sin flechas
  \tikzstyle{stateEdgePortionViolet} = [Plum,very thick,text= black];
  \tikzstyle{stateEdgeV} = [stateEdgePortionViolet];  % sin flechas
  \tikzstyle{stateEdgePortionGreen} = [ForestGreen,very thick,text= black];
  \tikzstyle{stateEdgeG} = [stateEdgePortionGreen];  % sin flechas
%  \tikzstyle{stateEdge} = [stateEdgePortion,->];  % con flechas
  \tikzstyle{edgeLabel} = [pos=0.5, text centered,very thick, font={\sffamily\small}];

\node[state, name=c1] {\underline{\textbf{Initial Data:}} Recurrence coefficients $a_n(t)$, $b_n(t)$};

\node[state, name=c2,below of=c1] {Discrete system \eqref{discrete system} and Toda system \eqref{toda system xn}};
\draw ($(c1.south) + (0em,0)$)
      edge[stateEdge] node[edgeLabel]{}
      ($(c2.north)$);

\node[state, name=c3,below of=c2] {Differential system \eqref{original system2}};
\draw ($(c2.south) + (0em,0)$)
      edge[stateEdge] node[edgeLabel]{}
      ($(c3.north)$);

\node[state, name=c4,below of=c3] {Points of indeterminacy \eqref{points}};
\draw ($(c3.south) + (0em,0)$)
      edge[stateEdge] node[edgeLabel]{}
      ($(c4.north)$);

\node[statepred, name=p1,below of=c4,left of=c4, xshift=-20em] {$P_1$};
\draw ($(c4.south) + (-2em,0)$)
      edge[stateEdge] node[edgeLabel]{}
      ($(p1.north)$);

\node[statepred, name=p5,below of=c4,left of=c4, xshift=-10em] {$P_5$};
\draw ($(c4.south) + (-1em,0)$)
      edge[stateEdge] node[edgeLabel]{}
      ($(p5.north)$);

\node[statepred, name=p2,below of=c4,left of=c4, xshift=-1em] {$P_2$};
\draw ($(c4.south) + (-0em,0)$)
      edge[stateEdge] node[edgeLabel]{}
      ($(p2.north)$);

\node[statepred, name=p3,below of=c4,right of=c4, xshift=-2em] {$P_3$};
\draw ($(c4.south) + (1em,0)$)
      edge[stateEdge] node[edgeLabel]{}
      ($(p3.north)$);

\node[statepred, name=p4,below of=c4,right of=c4, xshift=14em] {$P_4$};
\draw ($(c4.south) + (2em,0)$)
      edge[stateEdge] node[edgeLabel]{}
      ($(p4.north)$);

\node[state11, name=p17,below of=p1] {System $U_{1,1}V_{1,1}$ \eqref{systUV11original}};
\draw ($(p1.south) + (0em,0)$)
      edge[stateEdge] node[edgeLabel,xshift=-1em, yshift=0em]{\eqref{phi11}}
      ($(p17.north)$);
\draw ($(p5.south) + (-0.25em,0)$)
      edge[stateEdge] node[edgeLabel,xshift=-0.7em, yshift=0.5em]{$\otimes$}
      ($(p17.north)+ (+1em,0)$);

\node[state11, name=p26,below of=p5] {System $u_{5,4}v_{5,4}$ \eqref{systuv54op1}};
\draw ($(p5.south) + (0em,0)$)
      edge[stateEdge] node[edgeLabel,xshift=0.5em, yshift=0em]{$\oslash$}
      ($(p26.north)$);

\node[state11, name=p20,below of=p2] {System $U_{2,1}V_{2,1}$ \eqref{systUV21original}};
\draw ($(p2.south) + (0em,0)$)
      edge[stateEdge] node[edgeLabel,xshift=-1em, yshift=0em]{\eqref{phi21}}
      ($(p20.north)$);

\node[state11, name=p23,below of=p3] {System $U_{3,1}V_{3,1}$ \eqref{systUV31original}};
\draw ($(p3.south) + (0em,0)$)
      edge[stateEdge] node[edgeLabel,xshift=-1em, yshift=0em]{\eqref{phi31}}
      ($(p23.north)$);

\node[state11, name=p30,below of=p4] {System $u_{4,2}v_{4,2}$ \eqref{systuv42}};
\draw ($(p4.south) + (0em,0)$)
      edge[stateEdge] node[edgeLabel,xshift=3em, yshift=0em]{\eqref{phi41} and \eqref{phi42}}
      ($(p30.north)$);

\node[state11b, name=p52,below of=p17,right of=p17, xshift=-4em] {System $U_{1,2}V_{1,2}$ \eqref{systUV12}};
\draw ($(p17.south) + (1em,0)$)
      edge[stateEdgeG] node[edgeLabel,xshift=1.2em, yshift=0em]{\eqref{phi12}}
      ($(p52.north)$);

\node[state11b, name=p55,below of=p26,right of=p26, xshift=-4em] {System $u_{5,5}v_{5,5}$ \eqref{systuv552}};
\draw ($(p26.south) + (1em,0)$)
      edge[stateEdgeG] node[edgeLabel,xshift=1.2em, yshift=0em]{\eqref{phi55}}
      ($(p55.north)$);

\node[state11b, name=p53,below of=p20,right of=p20, xshift=-4em] {System $U_{2,2}V_{2,2}$ \eqref{systUV22}};
\draw ($(p20.south) + (1em,0)$)
      edge[stateEdgeG] node[edgeLabel,xshift=1.2em, yshift=0em]{\eqref{phi22}}
      ($(p53.north)$);

\node[state11b, name=p54,below of=p23,right of=p23, xshift=-4em] {System $U_{3,2}V_{3,2}$ \eqref{systUV32}};
\draw ($(p23.south) + (1em,0)$)
      edge[stateEdgeG] node[edgeLabel,xshift=1.2em, yshift=0em]{\eqref{phi32}}
      ($(p54.north)$);

\node[state11c, name=p43a,below of=p30,left of=p30, xshift=-0em,yshift=-2.4em] {System $u_{4,3^a}v_{4,3^a}$ \eqref{systuv431}\\$\shortparallel$\\System $U_{1,1}V_{1,1}$ \eqref{systUV11original}};
\draw ($(p30.south) + (-1em,0)$)
      edge[stateEdge] node[edgeLabel,xshift=-1.5em, yshift=0.5em]{\eqref{phi431}}
      ($(p43a.north)$);

\node[state11c, name=p43b,below of=p30, yshift=-2.4em] {System $u_{4,3^b}v_{4,3^b}$ \eqref{systuv432}\\$\shortparallel$\\System $U_{2,1}V_{2,1}$ \eqref{systUV21original}};
\draw ($(p30.south) + (-0em,0)$)
      edge[stateEdge] node[edgeLabel,xshift=-0.9em, yshift=0em]{\eqref{phi432}}
      ($(p43b.north)$);

\node[state11c, name=p43c,below of=p30,right of=p30, xshift=-0em,yshift=-2.4em] {System $u_{4,3^c}v_{4,3^c}$ \eqref{systuv433}\\$\shortparallel$\\System $U_{3,1}V_{3,1}$ \eqref{systUV31original}};
\draw ($(p30.south) + (+1em,0)$)
      edge[stateEdge] node[edgeLabel,xshift=1.5em, yshift=0.5em]{\eqref{phi433}}
      ($(p43c.north)$);

\node[state11d, name=p56,below of=p52] {Eq. $V_{1,2}$ \eqref{eq12}};
\draw ($(p52.south) + (0em,0)$)
      edge[stateEdgeG] node[edgeLabel,xshift=1.2em, yshift=0em]{\eqref{U12}}
      ($(p56.north)$);

\node[state11d, name=p43,left of=p56,xshift=2em] {Eq. $U_{1,1}$ \eqref{equationU11}};
\draw ($(p17.south) + (-1em,0)$)
      edge[stateEdgeV] node[edgeLabel,xshift=-1.2em, yshift=0em]{\eqref{V11}}
      ($(p43.north)$);

\node[state11d, name=p56b,below of=p55, yshift=-2.3em] {Eq. $v_{5,5}$ \eqref{systuv552}\\$\shortparallel$\\ Eq. $V_{1,2}$ \eqref{eq12}};
\draw ($(p55.south) + (0em,0)$)
      edge[stateEdgeG] node[edgeLabel,xshift=0.7em, yshift=0em]{\color{ForestGreen}{$\maltese$}}
      ($(p56b.north)$);

\node[state11d, name=p54b,left of=p56b,xshift=1.5em] {Eq. $v_{5,4}$ \eqref{equationv54}\\$\shortparallel$\\ Eq. $U_{1,1}$ \eqref{equationU11}};
\draw ($(p26.south) + (-1em,0)$)
      edge[stateEdgeV] node[edgeLabel,xshift=-0.7em, yshift=0em]{{\color{Plum}$\clubsuit$}}
      ($(p54b.north)$);

\node[state11d, name=p57,below of=p53] {Eq. $V_{2,2}$ \eqref{eq2.1}};
\draw ($(p53.south) + (0em,0)$)
      edge[stateEdgeG] node[edgeLabel,xshift=0.7em, yshift=0em]{{\color{ForestGreen}$\maltese$}}
      ($(p57.north)$);

\node[state11d, name=p44,left of=p57,xshift=2em] {Eq. $U_{2,1}$ \eqref{equationU21}};
\draw ($(p20.south) + (-1em,0)$)
      edge[stateEdgeV] node[edgeLabel,xshift=-0.7em, yshift=0em]{{\color{Plum}$\clubsuit$}}
      ($(p44.north)$);

\node[state11d, name=p58,below of=p54] {Eq. $V_{3,2}$ \eqref{eq31}};
\draw ($(p54.south) + (0em,0)$)
      edge[stateEdgeG] node[edgeLabel,xshift=0.7em, yshift=0em]{{\color{ForestGreen}$\maltese$}}
      ($(p58.north)$);

\node[state11d, name=p45,left of=p58,xshift=2em] {Eq. $U_{3,1}$ \eqref{equationU31}};
\draw ($(p23.south) + (-1em,0)$)
      edge[stateEdgeV] node[edgeLabel,xshift=-0.7em, yshift=0em]{{\color{Plum}$\clubsuit$}}
      ($(p45.north)$);

\node[state11e, name=f1,below of=p56,xshift=4em] {Painlev\'{e} V equation  with parameters \eqref{param1}};
\draw ($(p43.south) + (0em,0)$)
      edge[stateEdgeV] node[edgeLabel,xshift=-1em,yshift=-0.2em]{\eqref{transf moebius}}
      ($(f1.north) + (-2em,0)$);
\draw ($(p56.south) + (0em,0)$)
      edge[stateEdgeG] node[edgeLabel]{}
      ($(f1.north) + (-1em,0)$);
\draw ($(p54b.south) + (0em,0)$)
      edge[stateEdgeV] node[edgeLabel,xshift=-1em]{\eqref{transf moebius}}
      ($(f1.north) + (1em,0)$);
\draw ($(p56b.south) + (0em,0)$)
      edge[stateEdgeG] node[edgeLabel]{}
      ($(f1.north) + (2em,0)$);

\node[state11e, name=f2,right of=f1,xshift=1em] {Painlev\'{e} V equation  with parameters \eqref{param2}};
\draw ($(p44.south) + (0em,0)$)
      edge[stateEdgeV] node[edgeLabel,xshift=-1em,yshift=-0.2em]{\eqref{transf moebius}}
      ($(f2.north) + (-1em,0)$);
\draw ($(p57.south) + (0em,0)$)
      edge[stateEdgeG] node[edgeLabel]{}
      ($(f2.north) + (1em,0)$);

\node[state11e, name=f3,right of=f2,xshift=-1em] {Painlev\'{e} V equation  with parameters \eqref{param3}};
\draw ($(p45.south) + (0em,0)$)
      edge[stateEdgeV] node[edgeLabel,xshift=-1em,yshift=-0.2em]{\eqref{transf moebius}}
      ($(f3.north) + (-1em,0)$);
\draw ($(p58.south) + (0em,0)$)
      edge[stateEdgeG] node[edgeLabel]{}
      ($(f3.north) + (1em,0)$);

\node[state11c1, name=f4,below of=p43a] {Painlev\'{e} V equation  with parameters \eqref{param1}};
\draw ($(p43a.south) + (0em,0)$)
      edge[stateEdge] node[edgeLabel,xshift=-2em,yshift=0em]{Via  \eqref{systUV11original}}
      ($(f4.north) + (0em,0)$);

\node[state11c1, name=f5,right of=f4,xshift=-9em] {Painlev\'{e} V equation  with parameters \eqref{param2}};
\draw ($(p43b.south) + (0em,0)$)
      edge[stateEdge] node[edgeLabel,xshift=-2em,yshift=0em]{Via  \eqref{systUV21original}}
      ($(f5.north) + (0em,0)$);

\node[state11c1, name=f6,right of=f5,xshift=-9em] {Painlev\'{e} V equation  with parameters \eqref{param3}};
\draw ($(p43c.south) + (0em,0)$)
      edge[stateEdge] node[edgeLabel,xshift=-2em,yshift=0em]{Via  \eqref{systUV31original}}
      ($(f6.north) + (0em,0)$);

\end{tikzpicture}
}
\caption{Flow diagram.}
\label{diagrama}
\end{figure}

%%%%%%%%%%%%%%%%%%%%%%%%%%%%%%%%%%%%%%%%%%%%%%%%%%%%%%%%%%%%%%%%%%%%%%%%%%%%%%%%%%%%%%%

\begin{remark}
In the flow diagram presented in Figure~\ref{diagrama} we give a graphical representation of the regularisation procedure explained throughout Sections~\ref{regularisation} and~\ref{iterated regularisation}, as well as the ways we use to find a connection with the fifth Painlev\'{e} equation.  As it can be observed, starting with points $P_i$, $i=1,2,3$ we give two alternatives depending on the path chosen. The left path (in purple color) corresponds to apply only one iteration of the regularisation procedure while the right path (in green color) consists in applying a second iteration of the regularisation procedure which yields to polynomials systems (see Section~\ref{polynomial reg}). As expected, both paths allow finding connection to the fifth Painlev\'{e} equation with the same parameters in each case. We emphasise that both alternatives are of interest. On the one hand, the first path only requires a single iteration of the regularisation process followed by a M\"{o}ebius transformation. On the other hand, the second path involves two iterations of the regularisation procedure but no additional transformations, making this approach more algorithmic without the need of guessing any additional transformation. In this case, we do not observe a significant advantage between the two paths; however, we present both of them since in other situations involving more complicated initial systems, the second path might be more beneficial.

 For point $P_5$, we observe that, apart from the two paths analogous to those described for $P_i$, $i=1,2,3$ (paths in purple and green color), if we choose the left path in the first step, we arrive at the same path originated from $P_1$ (see Section~\ref{iterated regularisation}). Finally, for $P_4$ depending on the chosen path, we obtain the same systems as those derived from points $P_i$, $i=1,2,3$.

In each step (arrow) of the diagram the included labels make reference to the transformations which have been applied to obtain the following system or equation. In the case of $P_5$, we have not included all the transformations in the steps denoted by $\otimes$ and $\oslash$ to avoid overloading the diagram. The step denoted by $\otimes$ involves five iterations of the regularisation process (first one: transformations~(\ref{phi51b})--(\ref{phi54b}), second one: transformation (\ref{phi55}), third one: transformation (\ref{phi56}), fourth one: transformations (\ref{phi57b})--(\ref{phi510b}), fifth one: transformation (\ref{psi11})). On the other hand, the step denoted by $\oslash$ involves an iteration of the regularisation process consisting in transformations (\ref{phi51})--(\ref{phi54}).

The symbol {\color{Plum}$\clubsuit$} denotes transformations that have not been included in the paper to avoid an excessive number of transformations, but which can be derived following the same reasoning as for~(\ref{V11}). Similarly, the symbol {\color{ForestGreen}$\maltese$} denotes transformations that can be obtained reasoning in the same way as~(\ref{U12}).

Finally, we would like to remark that the arrows starting from equations~(\ref{eq12}),~(\ref{eq2.1}) and~(\ref{eq31}) do not involve any additinal transformations, since these equations coincide directly with the fifth Painlev\'{e} equation.

\end{remark}

\section{Hamiltonian expressions}
\label{hamiltonian}
Finally, we recall that the Painlev\'{e} equations can be written as Hamiltonian systems (see for instance~\cite[Chapter 32]{DMLF} and the references cited therein). Furthermore, we also have that~\eqref{systUV12},~\eqref{systUV22},~\eqref{systUV32},~\eqref{systuv552} can be expressed as Hamiltonian systems, that is, they satisfy the following expressions
\begin{equation*}
\begin{gathered}
\frac{dU_{i,2}}{dt}=\frac{\partial H_{i,2}}{\partial V_{i,2}}, \hspace{1 cm}
\frac{dV_{i,2}}{dt}=-\frac{\partial H_{i,2}}{\partial U_{i,2}},\\
\end{gathered}
\end{equation*}
for $i=1,2,3$, and
\begin{equation*}
\begin{gathered}
\frac{du_{5,5}}{dt}=\frac{\partial H_{5,5}}{\partial v_{5,5}}, \hspace{1 cm}
\frac{dv_{5,5}}{dt}=-\frac{\partial H_{5,5}}{\partial u_{5,5}},\\
\end{gathered}
\end{equation*}
 taking the functions
\begin{align*}
H_{1,2}=&\frac{V_{1,2} \left(N U_{1,2} \left(V_{1,2} (n-2 N-2)-N U_{1,2} \left(V_{1,2}-1\right){}^2+\alpha -n+2 N-t+2\right)\right)}{N t}\\
-&\frac{(N+1) (-n+N+1)}{Nt},\\
H_{2,2}=&\frac{N U_{2,2} \left(-V_{2,2} (-\alpha +n+2 N+t+2)+V_{2,2}^2 (n+N+1)-\alpha +N+1\right)-n (N+1) V_{2,2}}{N t}\\&-\frac{N^2 U_{2,2}^2 V_{2,2} \left(V_{2,2}-1\right){}^2}{Nt},\\
H_{3,2}=&\frac{-N U_{3,2} \left(V_{3,2} \left(V_{3,2} (\alpha -n+N+1)-\alpha +n-2 N+t-2\right)+N+1\right)+\alpha  V_{3,2} (n-N-1){}^2}{N t}\\
&-\frac{N^2 U_{3,2}^2 V_{3,2}
   \left(V_{3,2}-1\right)}{Nt},\\
H_{5,5}=&\frac{v_{5,5} \left(-u_{5,5} \left(v_{5,5} (n-2 N)+N u_{5,5} \left(v_{5,5}-1\right){}^2+\alpha -n+2 N-t\right)+n-N\right)+\alpha  u_{5,5}}{t},
\end{align*}
up to a function of $t$.\\

We also have that systems in other charts can be transformed to Hamiltonian polynomial systems. For instance, system~\eqref{systuv11original} in $(u_{1,1},v_{1,1})$ charts can be reduced to a polynomial system
\begin{equation}
\label{systuvoriginal}
\begin{gathered}
u_{1,2^b}'=\frac{-\alpha +n u_{1,2^b}-2 N u_{1,2^b}^2 v_{1,2^b}+2 N u_{1,2^b} v_{1,2^b}-t u_{1,2^b}^2+t u_{1,2^b}+\alpha  u_{1,2^b}}{t},\\
v_{1,2^b}'= \frac{-n N v_{1,2^b}+2 N^2 u_{1,2^b} v_{1,2^b}^2-N^2 v_{1,2^b}^2+2 N t u_{1,2^b} v_{1,2^b}-N t v_{1,2^b}-N t-\alpha
   N v_{1,2^b}-t}{N t},
\end{gathered}
\end{equation}
via the transformation
$$ u_{1,1}=-1+u_{1,2^b},\quad v_{1,1}=-\frac{1+N}{N}+u_{1,2^b}v_{1,2^b}.$$

Moreover, system~\eqref{systuvoriginal} can be expressed as
\begin{equation*}
\begin{gathered}
\frac{du_{1,2^b}}{dt}=\frac{\partial H_{1,2^b}}{\partial v_{1,2^b}}, \hspace{1 cm}
\frac{dv_{1,2^b}}{dt}=-\frac{\partial H_{1,2^b}}{\partial u_{1,2^b}},\\
\end{gathered}
\end{equation*}
for
$$H_{1,2^b}=\frac{v_{1,2^b} (u_{1,2^b} (\alpha +n-u_{1,2^b} (N v_{1,2^b}+t)+N v_{1,2^b}+t)-\alpha )}{t}+\frac{u_{1,2^b}}{N}+u_{1,2^b}.$$
Finally, we deduce the following statement for the original system~\eqref{original system2}.\\

\begin{theorem}The original system~\eqref{original system2} in $(q,p)$ coordinates can be expressed in the Hamiltonian form
\begin{equation*}
\begin{gathered}
\frac{1}{p+q}\frac{dq}{dt}=\frac{\partial H}{\partial p},\quad
\frac{1}{p+q}\frac{dp}{dt}=-\frac{\partial H}{\partial q},
\end{gathered}
\end{equation*}
for

$$H=\frac{N p q (\alpha +n-2 N-t-2)+N p^2 (n-N q)-q ((N+1) (-\alpha +N+1)+N t q)}{N t (p+q)}.$$

\end{theorem}

\section*{Discussion}

In this paper we analysed iterated regularisation for system~\eqref{original system2}, which helps to transform complicated systems with rational right--hand sides to simpler looking ones, even with polynomials right--hand sides. We also obtained that regularisation gives a nice decomposition of birational transformations. An intuitive idea for why this iterated polynomial regularisation process is effective is as follows. Starting from the initial system with rational right hand--side and by iterating the regularisation procedure, points of indeterminacy are progressively eliminated at each step. Eventually, this iterated process might transform the system into a system with polynomial right--hand side. Typically, it suffices to apply a simple M\"{o}bius transformation or a scaling to convert one of the final chart systems into a form equivalent to one of the canonical Painlev\'{e} equations. The main advantage of the iterated regularisation relies on being an algorithmic procedure (see Figure~\ref{diagrama}), avoiding the need of guessing the transformations required to establish connection with the Painlev\'{e} equation.

It is noticeable that alternative selections for the cascades of points may yield other simplified systems. Remarkably, the deduced equations from the polynomial systems coincide with Painlev\'{e} V equation.

\subsection*{Acknowledgements}

We sincerely thank the anonymous referees for their insightful comments and constructive
suggestions, which have greatly contributed to improve the clarity of the paper.

The author CRP would like to thank   Dr. hab Galina Filipuk for her  invitation to do a research stay at  the University of Warsaw.

\subsection*{Funding}

The authors are partially supported by the project PID2021--124472NB--I00 funded by \newline MCIN/AEI/10.13039/501100011033 and by ``ERDF A way of making Europe''. Furthermore, the work of the author CRP has been supported by PPIT--UAL, Junta de Andaluc\'{\i}a--ESF. Programme: 54.A. Aplication: 741. In addition, the work of the authors JFMM, JJMB and CRP  is partially supported by the project P\_FORT\_GRUPOS\_2023/72 from Plan Propio Investigaci\'{o}n de la Universidad de Almer\'{\i}a 2023, by the Research Group FQM--0229 and by the research centre CDTIME of Universidad de Almer\'{\i}a.

\subsection*{Competing Interests}
The authors have no relevant financial or non-financial interests to disclose.

\subsection*{Data Availability Statement }
Not applicable.

 \end{document}